\documentclass[12pt]{article}
\input{epsf.sty}
\usepackage{hyperref}
\usepackage{graphicx}
\usepackage{amsmath, amssymb}
\usepackage[arrow, matrix, curve]{xy}
\newtheorem {thm}{Theorem}[section]
\newtheorem {prop}[thm]{Proposition} 

\newtheorem {lem}[thm]{Lemma}
\newtheorem {cor}[thm]{Corollary}

\newtheorem {conj}{Conjecture}

\newcommand{\qed}{\nobreak \ifvmode \relax \else
      \ifdim\lastskip<1.5em \hskip-\lastskip
      \hskip1.5em plus0em minus0.5em \fi \nobreak
      \vrule height0.75em width0.5em depth0.25em\fi}

\def\Cox{\hfill \Box}

\def\N{{\Bbb N}}

\def\Z{{\Bbb Z}}
\def\R{{\Bbb R}}

\def\P{{\Bbb P}}

\def\E{{\Bbb E}}

\def\e{{\varepsilon}}

\def\eps{\varepsilon}

\def\ba{{\backslash}}

\def\D{\Delta}
\def\a{\alpha}

\def\ba{\setminus}
\def\b{\beta}

\def\d{\delta}

\def\e{\varepsilon}

\def\phi{\varphi}

\def\g{\gamma}

\def\l{\lambda}

\def\r{\rho}

\def\s{\sigma}

\def\o{\omega}

\def\D{\Delta}

\def\L{\Lambda}

\def\G{\Gamma}

\def\P{{\Phi}}

\def\T{\T}

\def\GG{{\cal G}}
\def\PP{{\cal P}}

\begin{document}
\title{Synchronization for discrete mean-field rotators 
}

\author{
Benedikt Jahnel
\footnote{ Ruhr-Universit\"at   Bochum, Fakult\"at f\"ur Mathematik, D44801 Bochum, Germany,
\newline
 \texttt{Benedikt.Jahnel@ruhr-uni-bochum.de}, 
\newline
\texttt{http://http://www.ruhr-uni-bochum.de/ffm/Lehrstuehle/Kuelske/jahnel.html }}
 \, and  Christof K\"ulske
\footnote{ Ruhr-Universit\"at   Bochum, Fakult\"at f\"ur Mathematik, D44801 Bochum, Germany,
\newline
\texttt{Christof.Kuelske@ruhr-uni-bochum.de}, 
\newline
\texttt{http://www.ruhr-uni-bochum.de/ffm/Lehrsttuehle/Kuelske/kuelske.html
/$\sim$kuelske/ }}\, 
\,  
}

\maketitle

\begin{abstract} 
We analyze a non-reversible mean-field jump dynamics for discrete $q$-valued rotators and show in particular 
that it exhibits synchronization. 
The dynamics is the mean-field analogue of the lattice dynamics investigated by the same authors in 
\cite{JK12} which   
provides an example of a non-ergodic interacting particle system on the basis of a mechanism 
suggested by Maes and Shlosman \cite{MS}.  

Based on the correspondence to an underlying model of continuous rotators via a discretization transformation
we show the existence of a locally attractive periodic orbit of rotating measures. 
We also discuss global attractivity, using a free energy as 
a Lyapunov function and the linearization of the ODE which describes typical behavior 
of the empirical distribution vector.




\end{abstract}



\vfill\eject

%
%
%
%
%
%
%
%
%
%
%
%

\section{Introduction}

Systems of interacting classical rotators ($S^1$-valued spins) on the sites of a lattice and also on different graphs 
have been a source of challenging and fruitful research in mathematical physics and probability. 
One likes to understand the nature of their translation-invariant phases (\cite{FSS76,BC99}), 
and the dependence on dimensionality (\cite{FrPf81}); one likes to understand the influence of different types of 
disorder, may it be destroying long-range order (\cite{AW90}) or even creating long-range order (\cite{C13});  
their dynamical properties, the difference that discretizations of the spin values make to the system 
(see the clock models in \cite{FS82}.) There is some similarity between rotators and massless models of real-valued unbounded fields (gradient fields), see  \cite{Fu05,EK08,CK12}. Roughly speaking the     
existence of ordered states for rotator models corresponds to existence of infinite-volume gradient states. 

There is usually much difference between the behavior of massless
models of continuous spins and models of 
discrete spins. 
The low energy excitations of the first are waves (see however the 
discrete symmetry breaking phenomenon of \cite{C13}), the  
excitations at very low temperatures of the latter can be described and controlled by contours (see \cite{B06}.)

There are however surprising situations when discrete models and continuous models behave the same: 
It is known that there can be even a continuum of extremal Gibbs measures 
for certain discrete-spin models (see \cite{FSS76}
for results in the nearest neighbor $q$-state clock 
model in an intermediate temperature regime.)  
 A route to create such a discrete system which is closely related but different 
 from the clock models with nearest neighbor interaction 
 goes as follows: Apply a sufficiently fine discretization transformation 
to the extremal Gibbs measures of an initial continuous-spin model in the regime where the initial system 
shows a continuous symmetry breaking.    Then show that 
the resulting uncountably many discretized measures are proper extremal Gibbs measures 
for a discrete interaction (see \cite{EKO11,JK12}.)
The model we are going to study here will also be of this type.  

There is another line of research leading to rotator models:   
Dynamical properties of rotator models from the rigorous and 
non-rigorous side have attracted a lot of interest from the statistical mechanics 
community and from the synchronisation community  (\cite{MS,BGP10,H00}.) 
Usually one studies a diffusive time-evolution of $S^1$- valued spins of mean-field type
which tends to synchronize the spins, where the mean-field nature is suggested by applications which come from systems of interacting neurons 
and collective motions of animal swarms.  Typically the dynamics is not reversible here. 
The first task one faces is to show (non-)existence of states describing collective synchronized 
motion, depending on parameter regimes. Next come questions about the 
approach of an initial state to these rotating states under time-evolution (\cite{BGP10,BGP12}), 
influence of the finite system size, 
and behavior 
at criticality (\cite{CD12}.) 

Our present research is motivated by a paper of Maes and Shlosman, \cite{MS}, 
about non-ergodicity in interacting particle systems (IPS). 
They conjectured 
that there could be non-ergodic behavior of a $q$-state IPS 
on the lattice in space dimensions $d\geq 3$ along the following mechanism involving rotating states.  
The system they considered was the $q$-state clock 
model with nearest neighbor scalarproduct interaction 
in an intermediate temperature regime where it is proved to 
have a continuity of extremal Gibbs states which can labelled by an angle.  
Then they proposed a  dynamics which 
should have the property to rotate the discrete spins according to local jump rules 
such that it possesses a periodic orbit consisting of these 
Gibbs states.Ê On the basis of  this heuristic idea of such a mechanism of rotating states,  
in a previous related work, \cite{JK12}, we considered a very special choice of quasilocal rates for  
a Markov jump process on the integer lattice in three or more spatial dimensions which provably 
shows this phenomenon. We were able to show that this IPS has 
a unique translation-invariant measure which is invariant under the dynamics but also possesses 
a non-trivial closed orbit of measures. Initialized at time zero according to a measure on this orbit 
the discrete spins perform synchronous rotations under the stochastic time evolution and 
don't settle in the time-invariant state. 
In particular we thereby constructed a lattice-translation invariant 
IPS which is non-ergodic in time. 
While such behavior was known to be possible for probabilistic cellular automata 
(infinite volume particle systems with simultaneous updating in discrete time), see \cite{CM11}, 
it was not known to occur for IPS (infinite volume particle systems in continuous time) and 
our example answers an old open question in IPS (Liggett question four of chapter one in \cite{Li85}.)  

There are open questions nonetheless in the lattice model. 
Of course it would be very interesting to see whether the periodic orbit of measures is attractive, 
what is the basin of attraction, what more can be said about the behavior of trajectories of time-evolved
measures, but this is open. 
We also don't know whether the original Maes-Shlosman conjecture is true and a simpler 
rotation dynamics with nearest neighbor interactions also behaves qualitatively the same 
in an intermediate temperature regime. 

In this paper let us therefore put ourselves to a mean-field situation and investigate whether 
we find analogies to the lattice and what more can be said now. 
This is interesting in itself since rotator models are naturally so often studied in a mean-field setting. 
What is a good version of a jump dynamics 
for discrete mean-field rotators implementing the Maes-Shlosman mechanism? 
Is there synchronisation for such a model as   
it is known to happen in the Kuramoto model (\cite{GPP12,DH96})? If yes, what can we say 
about attractivity of the orbit of rotating states?  Are there other attractors? 

Note that a very first naive attempt to define a discrete-spin mean-field dynamics showing 
synchonisation  does not work: 
the simple scalarproduct interaction $q$-state clock model does not have continuous symmetry breaking  at any $\b$. 
The model and its dynamics will rather 
appear as a discretization image of the continuous model on the level of measures. 
We consider the mean-field rotator model under equal-arc discretization into $q$ segments and  
define associated jump rates.   
Next we give criteria on the fineness of the discretization for existence and non-existence of the infinite-volume limit, 
and discuss a path large deviation principle (LDP) for empirical measures and the ODE for typical paths. 
 We prove that the discretization images of rotator Gibbs measures in the phase-transition region form a locally attractive limit cycle. Further we investigate local attractivity of the equidistribution and determine the non-attractive manifold. The question of global attractivity can be answered in the following way: Apart from measures with higher free energy than the equidistribution that get also trapped in the locally attractive manifold of the equidistribution, all measures are attracted by the limit cycle.

Summarizing, our mean-field results show many analogies to mean-field models of 
continuous rotators, they are in nice parallel to the behavior of the corresponding lattice system,  
but they go further since no stability result is known in the latter. 
It would be a challenge to see to what extend this parallel really holds.   

In the remainder of this introduction we present the construction and the main results without proofs.

\subsection{Model and Main Results}

We look at continuous-spin mean-field Gibbs measures in the finite volume $V_N=\{1,\dots, N\}$ 
which are the probability measures on the product space $(S^1)^N$ equipped with the product Borel sigma-algebra,  
defined by 
$$\mu_{\Phi,N}(d\s_{V_N})=\frac{\exp({-H_N(\s_{V_N})})\a^{\otimes N}(d\s_{V_N})}{\int_{(S^1)^N}\exp({-H_N(\bar\s_{V_N})})\a^{\otimes N}(d\bar\s_{V_N})}$$
where $\a$ is the Lebesgue measure on $S^1$. Here the energy function 
$$H_N(\s_{V_N})=N\Phi(L_N(\s_{V_N}))$$ depends on the spin configuration $\s_{V_N}=(\s_i)_{i\in V_N}$ 
only through the empirical distribution $L_N(\s_{V_N})=\frac{1}{N}\sum_{i=1}^N \d_{\s_i}$. 
Let us consider real-valued potentials $\Phi$ defined on the space of probability measures $\PP(S^1)$ on the sphere 
$S^1$ of two-body interaction type, 
$$\Phi(\nu)=\int\nu(ds_1)\int\nu(ds_2)V(s_1,s_2)$$
where $V$ is a symmetric pair-interaction function on $(S^1)^2$. We will refer to this model as the \textit{planar rotator model}. For the most 
part of the paper we will further specialize to the standard scalarproduct interaction with coupling strength $\b>0$ 
$$V(s_1,s_2)=-\frac{\b}{2}\langle e_{s_1},e_{s_2}\rangle$$ where $e_s$ is the unit vector pointing 
into the direction with angle $s$.
Recall as a standard fact that the distribution of the empirical measures $L_N$ 
under $\mu_{\Phi,N}$ obeys a LDP with rate $N$ 
and rate function given by the free energy 
$$\Psi(\nu )=\Phi(\nu)+S(\nu|\a)- \inf_{\tilde \nu}\bigl(\Phi(\tilde \nu)
+S(\tilde\nu|\a)\bigr)$$ where $S$ denotes the relative entropy. In the usual short notation let us write  
$$\mu_{\b,N}(L_N \approx \nu)\approx 
\exp \bigl(-N \Psi(\nu) \bigr).
$$
It is well known that there exist multiple minimizers of $\Psi$ in the scalarproduct model if and only if $\b>\b_c=2$ corresponding 
to a second-order phase transition in the inverse temperature at the critial value 
$2$ and a breaking of the $S^1$-symmetry. 
%
%
%
%
%
%
%

\subsubsection{Deterministic rotation, discretization and finite-volume Markovian dynamics for discretized systems}

For any real time $t$ we look at the joint rotation action $R_t:(S^1)^N\mapsto (S^1)^N$ given by the sitewise rotation 
of all spins, that is $(R_t\o_{V_{N}})_i=R_t \o_{i}$ where $R_t e_s=e_{(s+t)\text{mod}({2\pi})}$. 

Let  $\mu_N $ be a probability measure on $(S^1)^{N}$ which has a smooth 
Lebesque density relative to the product Lebesgue measure on  $(S^1)^{N}$. 
Denote the measure resulting from this deterministic rotation action $R_t$ by 
$\mu_{t, N}:=R_t\mu_{N}$.  

Next denote by T the local discretization map (local coarse-graining) with equal arcs of the sphere written as 
$[0,2 \pi)$ to the finite set $\{1,\dots, q \}Ê$, that is with  
 $S_k := [\frac{2 \pi}{q}(k-1),\frac{2 \pi}{q}k)$, $S^1=\bigcup_{k=1}^q S_k $ and  $T(s)=k$ if $s\in S_k$. 
 Extend this map to configurations in the product space by performing it sitewise. In particular we will consider images of 
measures under this discretization map $T$.

%
%
%

We will see that discretization after rotation of a continuous measure can be realized as a jump process.  
In order to define such a Markov jump process on the discrete-spin space $\{1,\dots, q\}^N$ 
we need some preparations. 
The following proposition describes the interplay between the discretization map $T$ and the deterministic 
rotation and is the starting point for the introduction of the dynamics we are going to consider.

\begin{prop}\label{Prop_Finite_Generator_Time_Dependent} There is a time-dependent linear generator $Q_{\mu_{N,t}}$ acting on discrete 
observables on the discrete $N$-particle state space,  $g:\{1,\dots, q\}^{N}\mapsto \R$,  such that 
an infinitesimal change of $T(\mu_{N,t})(g)= \int \mu_N(d \o)g( T R_t \o)$ can be written as 
\begin{equation}\label{Rotation_Finite_Volume}
\begin{split}
\lim_{\e \downarrow 0}\frac{1}{\e}\Bigl( T(\mu_{N,t+\e})(g)-T(\mu_{N,t})(g) \Bigr) 
=T(\mu_{N,t})(Q_{\mu_{N,t}} g). \cr
\end{split}
\end{equation}
This generator takes the form of a sum over single-site terms  
\begin{equation}\label{rotation_generator_finite_volume}
\begin{split}
Q_{\mu_{N,t}} g (\s'_{V_N})&:=\sum_{i=1}^N c_{\mu_{N,t}}(\s'_{V_N},\s_{V_N}'+1_i)(g(\s_{V_N}'+1_i)-g(\s'_{V_N}))\cr
\end{split}
\end{equation}
where $(\s_{V_N}'+1_i)_j =\s'_j + 1_{i=j}$ (modulo $q$). Here $c_{\mu_{N,t}}$ are certain time-dependent  
rates for increasing a coordinate by $1$ at single sites which have the feature to  
depend on time (only) through the measure $\mu_{N,t}$.\end{prop}

The generator $Q_{\mu_{N,t}}$ defines a Markov jump process (a continuous-time Markov chain) 
on the finite space $\{1,\dots, q\}^{N}$. There are only trajectories possible 
along which the variables $\s_i'$ 
increase their values by one unit along the circle of $q$ units according 
to the appropriate rates. An explicit expression for the rates 
in terms of the underlying measure can be found in formula \eqref{Rates_time_dependent}. 
The process we are going to study will be of this type. 

Specify to the case of a mean-field Gibbs measure $\mu_{\Phi,N}$ with a rotation-invariant interaction $\Phi$. 
Then $\mu_{\Phi,t,N}=\mu_{\Phi,N}$ stays constant under time-evolution 
and consequently the rates become time-independent. 
From  permutation invariance we see  that 
the resulting jump process obtained by mapping the trajectories 
of the paths $\s_{V_N}(t)$ to trajectories of the empirical distributions $L_N(\s_{V_N}(t))$ is a 
Markov process with generator which can be written in the form  
\begin{equation}\label{Rotation_Generator_Finite_Volume_Empirical}
\begin{split}
Q_N^{emp} f(\nu')&=N\sum_{k=1}^q\nu'(k)c_{N}^{emp}(k,\nu')( f(\nu'+\frac{1}{N}(\d_{k+1}-\d_{k}))-f(\nu')).\cr
\end{split}
\end{equation}
Here $f:\PP(\{1,\dots,q\})\mapsto \R$ is an observable on the simplex of 
$q$-dimensional probability vectors, 
$\d_k$ is the Dirac measure at $k\in \{1,\dots,q\}$
and $c_{N}^{emp}(k,\nu')$ are the resulting rates (given in \eqref{Finite_Volume_Rates_Conditional_Empirical}) describing the change of the empirical distribution at size $N$ 
when one particle changes its value from the state $k$ to $k+1$. 
 
 As a result of this construction of a Markovian dynamics we have the following corollary. 

\begin{cor}\label{Cor_Gibbs_Emp_Inv} Consider a mean-field Gibbs measure $\mu_{\Phi,N}$ for a rotation invariant potential $\Phi$. Then 
the stochastic dynamics on the space of empirical distributions $\PP(\{1,\dots,q\})$ 
with the above rates $c_{N}^{emp}(k,\nu')$ preserves the empirical distribution of the discretized mean-field Gibbs measure
$(T \mu_{\Phi,N})(L_N \in \,\cdot\,)   \in \PP(\{1,\dots,q\})$. \bigskip
\end{cor}

So far the construction of a mean-field dynamics for discrete rotators 
is largely in parallel to our construction of a dynamics for 
a non-ergodic IPS on $\Z^d$ as presented in \cite{JK12}. 

Our present aim for the mean-field setup is to understand large-$N$ properties, 
to understand mean-field analogues of rotating states and mean-field analogues of non-ergodicity. 
 We note that at finite $N$ of course we do not see a non-trivial closed orbit of measures. We will have to go to the limit $N\uparrow \infty$ to see reflections in the mean-field system of the non-ergodicity proved to occur for the IPS on the lattice. 
The picture one expects is the following: The empirical distribution (or profile) of a finite but very large particle system will become close in $\mathcal{O}(1)$ time to an empirical distribution (almost) on the periodic orbit. Then it will follow the orbit until a time large enough such that the finiteness of the system size will be felt. From that on it will not be sufficient to talk about a single profile anymore, rather more generally about a distribution of profiles, which, as time goes by, will mix over different angles along the orbit with equal probability. The relevant $N$-dependent mixing time we will not discuss in this paper. The control of closeness of the stochastic evolution up to finite times 
will be delivered by the path LDP which we are going to describe. Then we will analyze 
the typical behavior of the minimizing paths.   
While doing that we will be able to obtain additional information in mean field (which seem hard to get on the lattice) about stability of the periodic orbit under the dynamics.

\subsubsection{Infinite-volume limit of rates for fine enough discretizations}

To be able to understand the large-$N$ behavior we must look 
more closely to the rates $c_{N}^{emp}(k,\nu')$ and their large-$N$ limit. 
As it turns out, the existence and well-definedness is not completely automatic, 
but only holds if the discretization is sufficiently fine. This is an issue which is 
related to the appearance of non-Gibbsian measures under discretization transformations. 
On the constructive side we have the following result in our mean-field setup.

\begin{thm}\label{Thm_Rate_Convergence} For any smooth mean-field interaction potential $\Phi:\PP(S^1)\mapsto \R$ 
there is an integer $q(\Phi)$ such that for all $q\geq q(\Phi)$ 
the rates \eqref{Finite_Volume_Rates_Conditional_Empirical} have the infinite-volume limit 
\begin{equation*}\label{Rotation_Flow_Rates}
\begin{split}
c(k,\nu')&=\frac{\exp(-d\Phi_{\nu^{\nu'}})(\d_{\frac{2\pi}{q}k}-\nu^{\nu'})}
{\int_{S_k}\exp(-d\Phi_{\nu^{\nu'}}(\d_{\s}-\nu^{\nu'}))\a(d\s)}\cr
\end{split}
\end{equation*}
where the measure $\nu^{\nu'}$ is the unique
solution of the constrained free energy minimization problem 
$\nu \mapsto \Phi(\nu)+S(\nu|\a)$  in the set of $\nu\in \PP(S^1)$ 
with given discretization image $\nu'$, in other words in the set
$\{\nu\in \PP(S^1)|T(\nu)=\nu'\}$.
\end{thm}

Here $d\Phi_{\nu}(\d_{\s}-\nu)$ is the differential of the map $\Phi$ taken in the point $\nu\in \PP(S^1)$ applied 
to  the signed measure $\d_{\s}-\nu$ on $S^1$ with mass zero. 
It has the role of a mean field that a single spin feels when 
the empirical spin distribution in the system is $\nu$. 

The assumption of fine enough discretizations $q\geq q(\Phi)$ ensures that the
minimizer is unique and moreover Lipschitz continuous in total-variation distance 
as a function of $\nu'$ (see the proof of Lemma \ref{ODE_Solvable}.)  
The constrained minimizer $\nu^{\nu'}$ can be characterized as the unique solution 
of a typical mean-field consistency equation which reduces to a finite-dimensional 
equation in the case of the scalarproduct model. 
This uniqueness of the constrained free energy minimization  
is closely related to the notion of a mean-field Gibbs measure in terms of 
continuity of limiting conditional probabilities (see \cite{EKO11}.) 
The continuous spin value appearing in the definition of the rate to jump from $k$ to $k+1$ 
given by $\frac{2\pi}{q}k$ is the boundary between the segments of $S^1$ labelled by $k$ 
and by $k+1$. 
It is illuminating to compare the expression for the rates to the ones obtained 
for the non-ergodic IPS on the lattice from \cite{JK12} and observe 
the analogy. 

To get more concrete insight we specialize to the 
scalarproduct model where fineness criterion on discretization 
and form of rates are (more) explicit. We have the following proposition. 

\begin{thm}\label{Thm_Rate_Convergence_Rotator} Consider the standard scalarproduct model, 
let $\b>0$ be arbitrary (possibly in the phase-transition regime $\b>2$)
and $q$ be an integer large enough such that $\b\sin^2(\frac{\pi}{q})<1$.
Then the constrained free energy minimizer $\nu^{\nu'}$ is unique 
und the jump rates take the form 
\begin{equation}\label{Rotation_Flow_Rates}
\begin{split}
c(k,\nu')&=\frac{e^{\b\langle e_{\frac{2\pi}{q}k},M_\b(\nu')\rangle}}{\int_{S_k} e^{\b\langle e_\o,M_\b(\nu')\rangle}\a(d\o)}, \text{ for } k=0, \dots, q-1
\end{split}
\end{equation}
where $\nu'\mapsto M_\b(\nu'):=\int\nu^ {\nu'}(d\o)e_\o$ takes values in the two-dimensional unit disk. 
\end{thm}

The vector $M_\b(\nu')$ is the magnetization of the minimizing 
continuous-spin measure $\nu^{\nu'}$  which is constrained to $\nu'$. It is implicitly defined and 
can be computed from the solution of  
a mean-field fixed point equation.

The above criterion on the fineness of the discretization corresponds to 
the sufficient criterion for Gibbsianness of discretized lattice measures from \cite{KO}, \cite{EKO11}, \cite{JK12}. 
It is stronger than an application of the criterion for preservation of Gibbsianness under 
local transforms from \cite{KO08} would give (where however more general local transformations
were considered.)

We note that while some criterion on $q$ is necessary the present criterion 
is probably not sharp. Below we present an example where multiple constrained minimizers do actually occur (corresponding to non-Gibbssianness of the discretized model) which 
shows that large-$\beta$ asymptotics of the bound on $q$ is correct. 
The corresponding criterion is given in \eqref{Non_Gibbs_Criterion}. 

\subsubsection{Limiting dynamical system from path LDP as $N\uparrow \infty$}
 
 It is possible to formulate a path LDP for our dynamics. 
The infinite-volume limit of the rates enters into the
rate function. This rate function is a time-integral involving a 
Lagrangian density (see \eqref{Lagrangian_0}.) 
In the present introduction we restrict ourselves to formulate as a consequence 
the following (weak) law of large numbers (LLN) on the path level, for simplicity restricted to the planar rotor model. 

\begin{thm}\label{Path_LLN}
Let $\b\sin^2(\frac{\pi}{q})<1$, $T_f\in (0,\infty)$ be a finite time horizon. Let $(X_t)^N_{t\geq0}$
be the Markov jump process with generator $Q_N^{emp}$ started in an initial probability measure 
$\nu'_0$ on $\{1,\dots,q\}$. 
Then we have  $$(X_t)^N_{0\leq t\leq T_f}\xrightarrow[]{N\to\infty}(\phi(t,\nu'_0))_{0\leq t\leq T_f}$$ 
in the uniform topology on the pathspace, where the flow $\phi(t,\nu'_0)$ 
is given as a solution to 
the $(q-1)$-dimensional ordinary differential equation 
\begin{equation}\label{Rotation_Flow}
\begin{split}
\frac{d}{dt}\phi(t,\nu'_0)&=F(\phi(t,\nu'_0))
\end{split}
\end{equation}
with initial condition $\phi(0,\nu'_0)=\nu'_0$,  
for the vector field $F(\nu')=(F(\nu')(k))_{k=1,\dots,q}$ acting on $\PP(\{1,\dots,q\})$ with components
\begin{equation}\label{Rotation_Flow_RHS}
\begin{split}
F(\nu')(k)&=c(k-1,\nu')\nu'(k-1)-c(k,\nu')\nu'(k), \qquad k=1,\dots, q.
\end{split}
\end{equation}
 \end{thm}
While the LLN could also be obtained differently (and maybe more 
easily) the LDP from which this result follows is of independent 
interest of course.  It provides an interesting link with Lagrangian dynamics. 
 Its proof uses the Feng-Kurtz scheme (see \cite{FK06}.)   

The dynamical system with vector field $F$ introduced above provides the mean-field analogue 
in the large-$N$ limit of the non-ergodic IPS from \cite{JK12}.  
So one expects that it should reflect the non-ergodic lattice behavior 
(based on the rotation of states) by showing a closed orbit and we will see that this is really the case.

\subsubsection{Properties of the flow: Closed orbit and equivariance property of the discretization map}

Now we will come to the discussion of the analogue of the breaking of ergodicity in the IPS in \cite{JK12} occuring on the level of the infinite-volume limit of the mean-field system. 
Denote the continuous-spin free energy minimizers (infinite-volume Gibbs measures on empirical magnetization) by 
$$\GG(\Phi):=\text{argmin}(\nu \mapsto \Psi(\nu)).$$ 
Denote the discrete-spin free energy minimizers by the measures 
$$\GG':=\text{argmin}(\nu' \mapsto \Psi'(\nu'))$$
where the discrete-spin free energy function $\Psi'$
\begin{equation*}\label{Free_Energy}
\begin{split}
\Psi'(\nu')&:=\Psi(\nu^{\nu'})
\end{split}
\end{equation*}
is defined by the constrained minimization.
 
The vector field $F$ has the property that deterministic rotation 
of free energy minimizers in $\PP(S^1)$ is reproduced by the flow of free energy minimizers in $\PP(\{1,\dots, q\})$.
In the phase-transition regime of the planar rotor model 
the continuous-spin free energy minimizers in $\PP(S^1)$ 
can be labelled by the angle of the magnetization values. 
Hence the vector field $F$ has a closed orbit. We can summarize the 
interplay between discretization, deterministic rotation of continuous measures 
and evolution according to the flow $(\phi_t)_{t\geq0}$ of 
the ODE for discrete measures in the following picture. 

\begin{thm}\label{Diagram_Commutating}

The following diagram is commutating
$$
\begin{xy}
  \xymatrix{
      \PP(S^1)\supset&\GG(\Phi) \ar[rrr]^{\nu\mapsto R_t\nu} \ar@/_0,5cm/[d]_T &  &  &\GG(\Phi) \ar[d]^T \\
      \PP(\{1,\dots,q\})\supset&\GG' \ar[rrr]_{\nu'\mapsto \phi_t\nu'} \ar@/_0,5cm/[u]_{\nu'\mapsto \nu^{\nu'}}        &   & &\GG' 
  }
\end{xy}
$$
\end{thm}

This picture is in perfect analogy to the behavior of the IPS from \cite{JK12}. 
(Let us  point out that the generator from \cite{JK12} is more involved 
since it contains another part corresponding 
to a Glauber dynamics. This part was added for reasons which are not present 
in the mean-field setup. 
It will not be treated here.)

\subsubsection{Properties of the flow: Attractivity of the closed orbit }

For the following we restrict to the standard scalarproduct model and we assume 
that we are in the regime $\b>2$ where a non-trivial closed orbit exists.   
We want to understand the dynamics in the infinite-volume limit. In our present mean-field setup this boils down to a discussion of the finite-dimensional ODE, so we are left at this stage with a purely 
analytical question. Note that our ODE for discrete rotators 
parallels a non-linear PDE for the continuous rotators with all its 
intricacies (see \cite{BGP12}.) Having the benefit of finite dimensions however we have to 
deal with the additional difficulty that in our case the r.h.s is only implicitly defined. 

As our dynamics is non-reversible it is not clear a priori what 
the behavior of the free energy $\Psi'$ for the discrete system 
will be under time evolution. 
However, since we already know that the ODE has as a periodic orbit, namely the set of discretization images of continuous free energy minimizers, we might hope that the free energy $\Psi'$ will work as a Lyapunov function. As it turns out this is the case.

\begin{prop}\label{Lyapunov_Lemma}
Under the flow $\phi(t,\nu')$ the discrete-spin free energy $\Psi'$ is non-increasing, 
$\frac{d}{dt}\bigl |_{t=0}\Psi'(\phi(t,\nu'))\leq 0$,  for all $\nu'\in \PP(\{1,\dots,q\})$. 
The free energy does not change, 
$\frac{d}{dt}\bigl |_{t=0}\Psi'(\phi(t,\nu'))= 0$, 
 if and only if  $\nu'\in\GG'$ or $\nu'=\frac{1}{q}\sum_{k=1}^q \d_{k}$.

%

%

\end{prop}

The proof is not as obvious as one would hope for and uses
change of variables to new variables after which certain convexity properties 
can be used.   
This seems to be particular to the standard scalarproduct model.  
As a corollary we have the attractivity of the periodic orbit formulated as follows. 
\begin{thm}\label{Attract_Result} 
For any starting measure $\nu'\in \PP(\{1,\dots,q\})$ 
with free energy \\
$\Psi'(\nu')<\Psi'(\frac{1}{q}\sum_{k=1}^q \d_{k})$ the trajectory $\phi(t,\nu')$ enters
any open neighborhood around the periodic orbit $\GG'$ after finite time $t$. 

\end{thm}

\subsubsection{Properties of the flow: Stability analysis at the equidistribution}

\begin{figure}[h]
\begin{center}
\includegraphics[width=14cm]{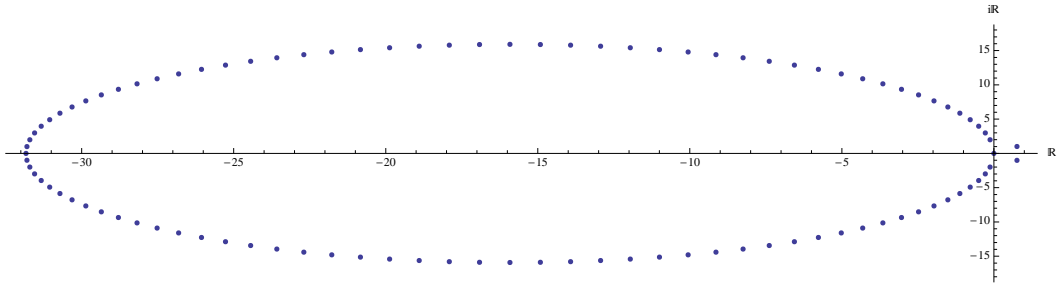}
\end{center}
\caption{\scriptsize{Spectrum of $\tilde M$ for $q=100, \b=50$. The two positive real-part eigenvalues are clearly visible.}}
\label{Spectrum}
\end{figure}

For the case of initial conditions 
$\nu'$ with 
free energy   $\Psi'(\nu')\geq\Psi'(\frac{1}{q}\sum_{k=1}^q \d_{k})$ 
we only know from 
the previous reasoning that 
the trajectories enter any open neighborhood around periodic orbit {\em and} equidistribution 
after finite time.  So we are interested in the stability of the dynamics locally 
around the equidistribution. Computing the linearization of the r.h.s of the ODE from its defining  
fixed point equation and using discrete Fourier transform  
we derive explicit expressions for its eigenvalues (see Lemma \ref{Eigenvalues_Lemma} and figure~\ref{Spectrum}.) We see that the linearized dynamics  rotates and exponentially 
suppresses the discrete Fourier-modes of the empirical measure 
except the lowest one which is expanded. In particular we have the following result which is in analogy to the behavior of 
the continuous model of \cite{GPPP12}.   
\begin{thm}\label{Attract_Equi} 
In the relevant parameter regimes (i.e. outside the set of parameters satisfying criterion \eqref{Non_Gibbs_Criterion}), the equidistribution is locally not purely attractive. The $2$-dimensional non-attractive manifold is given by
\begin{equation*}\label{Attractive_Manifold}
 \begin{split}
 \Bigl\{\nu'\in \PP(\{1,\dots,q\})\Bigl |\sum_{k=1}^q\nu'(k)e^{i\frac{2\pi}{q}lk}=0\text{ for all }l\in\{2,\dots,q-2\}\Bigr\}.
\end{split}
\end{equation*}
\end{thm}


\begin{figure}[h]
\begin{center}
\includegraphics[width=14cm]{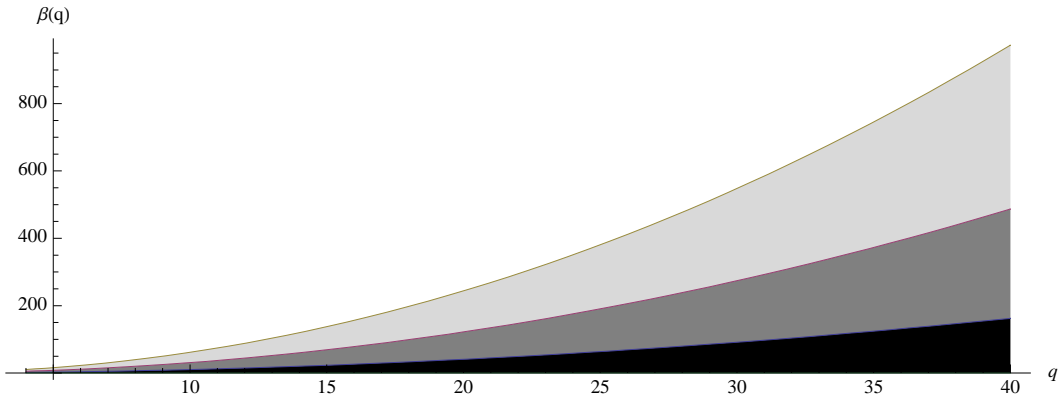}
\end{center}
\caption{\scriptsize{The black area shows $(\b,q)$-regimes where uniqueness of constrained free energy minimizers is guaranteed by the criterion given in Theorem \ref{Thm_Rate_Convergence_Rotator}, in other words our construction certainly works. The light gray and white areas show $(\b,q)$-regimes where the complemetary criterion \eqref{Non_Gibbs_Criterion} holds, in other words our limiting dynamics can not be defined. In the intermediate dark gray area we do not know whether our dynamics is well defined. However only in the white area the equidistribution is purely attractive, in the relevant 
$(\b,q)$-regimes we have non-attractivity. All analysis is done for $\b>2$ since we only work in the phase-transition region.}}
\label{Plot_Parameter_Regimes}
\end{figure}

\subsubsection{Outline of the paper}

In section \eqref{Rotation_Dynamics} subsection \ref{Finite-volume rotation dynamics} we consider the rotation dynamics first in the finite volume as an IPS on the level of spins. We prove Proposition \ref{Prop_Finite_Generator_Time_Dependent}. After specializing to the case where the finite-volume dynamics leaves the Gibbs measure invariant, we lift the dynamics to the level of empirical distributions and prove Corollary \ref{Cor_Gibbs_Emp_Inv}. In subsection \ref{The infinite-volume rotation rates: Existence and non-existence} we prove Theorem \ref{Thm_Rate_Convergence} on the convergence of the rates in the thermodynamic limit. We prove a useful lemma about uniqueness of constrained free energy minimizers for fine enough discretization, still for a more general interaction potential. Here we follow an adaptation of arguments presented in \cite{EKO11,G11,KO08}. The proof of Theorem \ref{Thm_Rate_Convergence_Rotator} uses the special structure of the standard scalarproduct interaction to derive a tangible criterion for the fineness of discretization implying uniqueness 
of constrained minimizers which are needed for the existence of limiting rates for the dynamics. 
In \eqref{Non_Gibbs_Criterion} we present a complementary criterion on 
the coarseness of the discretization ensuring non-uniqueness of constrained minimizers.  
In subsection \ref{The infinite-volume rotation dynamics} we prove global existence of solutions of the infinite-volume dynamics via Lipschitz continuity of the r.h.s. Further we prove Theorem \ref{Path_LLN} employing a LDP on the level of paths.

\bigskip
Section \ref{Properties of the flow} subsection \ref{Properties of the flow: Closed orbit and equivariance property of the discretization map} contains the proof of the equivariance property indicated in the diagram of Theorem \ref{Diagram_Commutating}. 
In subsection \ref{Attractivity via the free energy} we derive the time-derivative of the free energy and prove Proposition \ref{Lyapunov_Lemma}. 
As a consequence we obtain stability of the periodic orbit formulated in Theorem \ref{Attract_Result}.
Subsection \ref{Local stability analysis of the equidistribution via linearization} is devoted to the local stability analysis at the equidistribution and the proof of Theorem \ref{Attract_Equi}. 

\bigskip
\textbf{Acknowledgement: }  
We thank A. Abbondandolo, A. van Enter,  G. Giacomin,  G. Knieper, 
C. Maes, A. Opoku, F. Redig, and W. Ruszel for stimulating discussions.  
This work is supported by the Sonderforschungsbereich SFB $|$ TR12-Symmetries and Universality in Mesoscopic Systems.

\section{Rotation dynamics}\label{Rotation_Dynamics}


\subsection{Finite-volume rotation dynamics}\label{Finite-volume rotation dynamics}

We consider the time-dependent generator \eqref{rotation_generator_finite_volume}
acting on discrete observables on the discrete $N$-particle state space
and $\mu_{N,t}(d\s_{V_N})=\r_{N,t}(\s_{V_N})\a^{\otimes N}(d\s_{V_N})$ where $\a$ denotes the Lebesgue measure on $S^1$ and the density $\r_{N,t}$ is supposed to be continuous. The time-dependent rates are given by 
\begin{equation}\label{Rates_time_dependent}
\begin{split}
c_{\mu_{N,t}}(\s'_{V_N},\s'_{V_N}+1_i)&:=\frac{\int_{T^{-1}(\s'_{V_N\ba i})}\r_{N,t}(\frac{2\pi}{q}\s'_i,\s_{V_N\ba i})\a^{\otimes N\ba i}(d\s_{V_N\ba i})}{\int_{T^{-1}(\s'_{V_N})}\r_{N,t}(\s_{V_N})\a^{\otimes N}(d\s_{V_N})}.\cr
\end{split}
\end{equation}


\bigskip
\textbf{Proof of Proposition \ref{Prop_Finite_Generator_Time_Dependent}: }
It suffices to check \eqref{Rotation_Finite_Volume} for $g=1_{\s'_{V_N}}$. We have
\begin{equation*}\label{Key_Derivative}
\begin{split}
&\lim_{\e \downarrow 0}\frac{1}{\e}\Bigl( T(\mu_{N,t+\e})(g)-T(\mu_{N,t})(g) \Bigr)\cr
&=\lim_{\e \downarrow 0}\frac{1}{\e}\Bigl( \frac{\int_{T^{-1}(\s'_{V_N})}\r_{N,t+\e}(\s_{V_N})\a^{\otimes N}(d\s_{V_N})}{\int\r_{N,t+\e}(\s_{V_N})\a^{\otimes N}(d\s_{V_N})}-\frac{\int_{T^{-1}(\s'_{V_N})}\r_{N,t}(\s_{V_N})\a^{\otimes N}(d\s_{V_N})}{\int\r_{N,t}(\s_{V_N})\a^{\otimes N}(d\s_{V_N})} \Bigr)\cr
&=\frac{1}{Z}\sum_{i=1}^{N}\lim_{\e \downarrow 0}\frac{1}{\e}\Bigl( \int_{T^{-1}(\s'_{V_N\ba i})}\bigl(\int_{\frac{2\pi}{q}(\s'_i-1)-\e}^{\frac{2\pi}{q}(\s'_i-1)}\r_{N,t}(\s_{V_N})-\int^{\frac{2\pi}{q}\s'_i}_{\frac{2\pi}{q}\s'_i-\e}\r_{N,t}(\s_{V_N})\bigr)\a^{\otimes N}(d\s_{V_N})\Bigr)\cr
&=\sum_{i=1}^{N}\Bigl(c_{\mu_{N,t}}(\s'_{V_N}-1_i,\s'_{V_N})T(\mu_{N,t})(\s'_{V_N}-1_i)-c_{\mu_{N,t}}(\s'_{V_N},\s'_{V_N}+1_i)T(\mu_{N,t})(\s'_{V_N})\Bigr)\cr
&=T(\mu_{N,t})(Q_{\mu_{N,t}} g) \cr
\end{split}
\end{equation*}
$\Cox$

Plugging in for $\r$ the Gibbs density for a rotation-invariant potential, the rates take the time-independent form
\begin{equation}\label{Gibbs_Rates}
\begin{split}
c_{\mu_{\Phi,N}}(\s'_{V_N},\s'_{V_N}+1_i)&=\frac{\int_{T^{-1}(\s'_{V_N\ba i})}e^{-N\Phi(\frac{1}{N}\d_{\frac{2\pi}{q}\s'_i}+\frac{N-1}{N}L_{N-1}(\s_{V_N\ba i}))}\a^{\otimes N\ba i}(d\s_{V_N\ba i})}{\int_{T^{-1}(\s'_{V_N})}e^{-N\Phi(L_{N}(\s_{V_N}))}\a^{\otimes N}(d\s_{V_N})}\cr
\end{split}
\end{equation}
and $T(\mu_{\Phi,N})(Q_{\mu_{\Phi,N}}g)=0$ for all discrete observables $g$. Hence $T(\mu_{\Phi,N})$ is invariant under $Q_{\mu_{\Phi,N}}$.
Notice one can rewrite the rates as
\begin{equation*}\label{Finite_Volume_Rates_Conditional}
\begin{split}
c_{\mu_{\Phi,N}}&(\s'_{V_N},\s'_{V_N}+1_i)=\cr
&\frac{\mu_{\Phi,N-1}[\s'_{V_{N}\ba i}](e^{-N\Phi(\frac{1}{N}\d_{\frac{2\pi}{q}\s'_i}+\frac{N-1}{N}L_{N-1}(\cdot))+(N-1)\Phi(L_{N-1}(\cdot))})}{\mu_{\Phi,N-1}[\s'_{V_{N}\ba i}](\int_{T^{-1}(\s'_i)}e^{-N\Phi(\frac{1}{N}\d_{\s_i}+\frac{N-1}{N}L_{N-1}(\cdot))+(N-1)\Phi(L_{N-1}(\cdot))}\a(d\s_i))}
\end{split}
\end{equation*}
where $\mu_{\Phi,N-1}[\s'_{V_{N}\ba i}]$ stands for the Gibbs measure conditioned to the set $T^{-1}(\s'_{V_{N}\ba i})$. In fact only empirical distributions of the coarse-grained spin variables $L_{N-1}(\s'_{V_{N}\ba i})$ come into play. Thus by writing $\nu'\in \PP(\{1,\dots,q\})$ for a possible empirical measure $L_{N-1}$ we can again re-express the rates as
\begin{equation}\label{Finite_Volume_Rates_Conditional_Empirical}
\begin{split}
c_N^{emp}(k,\nu')&=\frac{\mu_{\Phi,N-1}^{emp}[\nu'](e^{-N\Phi(\frac{1}{N}\d_{\frac{2\pi}{q}
k}+\frac{N-1}{N}L_{N-1}(\cdot))+(N-1)\Phi(L_{N-1}(\cdot))})}{\mu_{\Phi,N-1}^{emp}[\nu'](\int_{S_k}e^{-N\Phi(\frac{1}{N}\d_{\s}+\frac{N-1}{N}L_{N-1}(\cdot))+(N-1)\Phi(L_{N-1}(\cdot))}\a(d\s))} 
\end{split}
\end{equation}
where we now dropped the indication for the Gibbs measure. Notice, this expression makes sense even when $\nu'$ is not an empirical distribution.

We can now lift the whole process to the level of empirical distributions. The resulting generator is given in \eqref{Rotation_Generator_Finite_Volume_Empirical}.

\bigskip
\textbf{Proof of Corollary \ref{Cor_Gibbs_Emp_Inv}: }
We have to check $T(\mu_{\Phi,N})((Q_N^{emp}f)\circ L_N)=0$ for all bounded measurable functions $f:\PP(\{1,\dots,q\})\mapsto\R$.
\begin{equation*}\label{Lift}
\begin{split}
&T(\mu_{\Phi,N})((Q_N^{emp}f)\circ L_N)
=\sum_{\s'_{V_N}}T(\mu_{\Phi,N})(\s'_{V_N})\sum_{i=1}^N\sum_{k=1}^q\d_{\s'_i=k}c_N^{emp}(k,L_N(\s'_{V_N}))\times\cr
&\hspace{6cm}\bigl(f(L_N(\s'_{V_N})+\frac{1}{N}(\d_{k+1}-\d_{k}))-f(L_N(\s'_{V_N}))\bigr)\cr
&=\sum_{\s'_{V_N}}T(\mu_{\Phi,N})(\s'_{V_N})\sum_{i=1}^Nc_{\mu_{\Phi,N}}(\s'_{V_N},\s'_{V_N}+1_i)\bigl(f(L_N(\s'_{V_N}+1_i))-f(L_N(\s'_{V_N}))\bigr)\cr
&=T(\mu_{\Phi,N})(Q_{\mu_{\Phi, N}}(f\circ L_N)).\cr
\end{split}
\end{equation*}
But $T(\mu_{\Phi,N})(Q_{\mu_{\Phi, N}}(f\circ L_N))=0$ since $T(\mu_{\Phi,N})$ is invariant for $Q_{\mu_{\Phi, N}}$.
$\Cox$

\subsection{Infinite-volume rates: Existence and non-existence}\label{The infinite-volume rotation rates: Existence and non-existence}
Let us prepare the proof of Theorem \ref{Thm_Rate_Convergence} by the following lemma.
\begin{lem}\label{Fine_Descretization_General}
For any differentiable mean-field interaction potential $\Phi:\PP(S^1)\mapsto \R$ with
\begin{equation*}\label{Lipschitz_Derivative of potential}
\begin{split}
\sup_{s,t\in S_k}|d_{\mu}\Phi(\d_s-\d_t)-d_{\tilde\mu}\Phi(\d_s-\d_t)|\leq C(q)\Vert\tilde\mu-\mu\Vert
\end{split}
\end{equation*}
where $C(q)\downarrow0$ for $q\uparrow\infty$ monotonically, there is an integer $q(\Phi)$ such that for all $q\geq q(\Phi)$  the free energy minimization problem $\nu \mapsto \Phi(\nu)+S(\nu|\a)$ has a unique solution in the set $\{\nu\in \PP(S^1)|T(\nu)=\nu'\}$ for any $\nu'\in\PP(\{1,\dots,q\})$. 
\end{lem}
We call this solution $\nu^{\nu'}$. The proof follows a line of arguments given in \cite{EKO11} in the lattice situation.

\bigskip
\textbf{Proof of Lemma \ref{Fine_Descretization_General}: }
Let $\mu$ be a solution of the constrained free energy minimization problem $\nu \mapsto \Phi(\nu)+S(\nu|\a)$ with $T(\mu)=\nu'$ and  $\tilde\mu$ be a solution of the constrained free energy minimization problem $\nu \mapsto \tilde\Phi(\nu)+S(\nu|\a)$ with $\tilde\Phi$ being another continuously differentiable mean-field interaction potential and $T(\tilde\mu)=\nu'$. Using Lagrange multipliers to characterize the constrained extremal points of the free energy we find $\mu$ and $\tilde\mu$ must have the form
\begin{equation}\label{MF_DLR_General}
\begin{split}
\mu(ds|S_k)&=\frac{1_{S_k}\exp(-d_\mu\Phi(\d_s-\mu))}{\int_{S_k}\exp(-d_\mu\Phi(\d_{\bar s}-\mu))\a(d\bar s)}\a(ds)=:\g_k(ds|\mu)\cr
\tilde\mu(ds|S_k)&=\frac{1_{S_k}\exp(-d_{\tilde\mu}\tilde\Phi(\d_s-\tilde\mu))}{\int_{S_k}\exp(-d_{\tilde\mu}\tilde\Phi(\d_{\bar s}-\tilde\mu))\a(d\bar s)}\a(ds)=:\tilde\g_k(ds|\tilde\mu).
\end{split}
\end{equation}
Let us estimate for a bounded measurable function $f$
\begin{equation}\label{MF_EKO}
\begin{split}
|\mu(f|S_k)-\tilde\mu_k(f|S_k)|&\leq|\g_k(f|\mu)-\g_k(f|\tilde\mu)|+|\g_k(f|\tilde\mu)-\tilde\g_k(f|\tilde\mu)|.\cr
\end{split}
\end{equation}
With $\Vert\a_1-\a_2\Vert:=\max_{f \text{ bounded, measurable}}|\a_1(f)-\a_2(f)|/\d(f)$ denoting the total-variation distance of probability measures where $\d(f):=\sup_{x,y}|f(x)-f(y)|$ is the variation of a bounded function we have
\begin{equation*}\label{MF_EKO_Second}
\begin{split}
|\g_k(f|\tilde\mu)-\tilde\g_k(f|\tilde\mu)|\leq\d(f)\Vert\g_k(\cdot|\tilde\mu)-\tilde\g_k(\cdot|\tilde\mu)\Vert=:\d(f)b(\tilde\mu).\cr
\end{split}
\end{equation*}
For the first term in \eqref{MF_EKO} we similary write 
\begin{equation*}\label{MF_EKO_Second}
\begin{split}
|\g_k(f|\mu)-\g_k(f|\tilde\mu)|\leq\d(f)\Vert\g_k(\cdot|\mu)-\g_k(\cdot|\tilde\mu)\Vert.\cr
\end{split}
\end{equation*}
Now let $u_1(s):=d_\mu\Phi(\d_s-\mu)$,  $u_0(s):=d_{\tilde\mu}\Phi(\d_s-\tilde\mu))$, $v:=u_1-u_0$ and $u_t:=u_0+tv$. Define $h_t^k:=\exp(u_t)1_{S_k}/\a(\exp(u_t)1_{S_k})$ and $\l_t^k(ds):=h_t^k(s)\a(ds)$. Then we have
\begin{equation}\label{MF_EKO_key_1}
\begin{split}
2\Vert\g_k(\cdot|\mu)-\g_k(\cdot|\tilde\mu)\Vert&=2\Vert\l_1^k-\l_0^k\Vert\leq\int_0^1dt\a(|\frac{d}{dt}h_t^k|)\cr
&=\int_0^1dt\l_t^k(|v-\l_t^k(v)|)\cr
&\leq\int_0^1dt\int\l_t^k(dx)\int\l_t^k(dy)|v(x)-v(y)|\cr
&=\int_0^1dt\int v(\l_t^k)(dx)\int v(\l_t^k)(dy)|x-y|\cr
&\leq\sup_\l\int_{-r}^r\l(dx)\int_{-r}^r\l(dy)|x-y|\cr
\end{split}
\end{equation}
where the supremum is over all probability measures on the interval $[-r,r]$ with 
$2r:=\sup_{s,t\in S_k}|d_{\mu}\Phi(\d_s-\d_t)-d_{\tilde\mu}\Phi(\d_s-\d_t)|$. By assumption
we have
\begin{equation}\label{MF_EKO_Lipischitz_Absch}
\begin{split}
2r\leq C(q)\Vert\mu-\tilde\mu\Vert\leq C(q)\sup_{l\in\{1,\dots,q\}}\Vert\mu(\cdot|S_l)-\tilde\mu(\cdot|S_l)\Vert
\end{split}
\end{equation}
with $C(q)\downarrow0$ for $q\uparrow\infty$ monotonically. Using the fact, that for all probability measures $p$ on $[-r, r]$ we have $\int p(dx)\int p(dy )|x-y|\leq r$ and \eqref{MF_EKO_Lipischitz_Absch} we can thus find $q(\Phi)$ such that
\begin{equation*}\label{MF_EKO_key_2}
\begin{split}
\Vert\g_k(\cdot|\mu)-\g_k(\cdot|\tilde\mu)\Vert\leq C(q(\Phi))\sup_{l\in\{1,\dots,q\}}\Vert\mu(\cdot|S_l)-\tilde\mu(\cdot|S_l)\Vert
\end{split}
\end{equation*}  
with $C(q(\Phi))<1$. Hence for all $q\geq q(\Phi)$
\begin{equation*}\label{MF_EKO_2}
\begin{split}
|\mu(f|S_k)-\tilde\mu(f|S_k)|&\leq\d(f)(C(q)\sup_{l\in\{1,\dots,q\}}\Vert\mu(\cdot|S_l)-\tilde\mu(\cdot|S_l)\Vert+b(\tilde\mu)).\cr
\end{split}
\end{equation*}
Taking the supremum over $f$ and over $k$ we have 
\begin{equation*}\label{MF_EKO_3}
\begin{split}
\sup_{k\in\{1,\dots,q\}}\Vert\mu(\cdot|S_k)-\tilde\mu(\cdot|S_k)\Vert&\leq\frac{1}{1-C(q)}b(\tilde\mu).\cr
\end{split}
\end{equation*}
Now for $\tilde\Phi=\Phi$ of course $b(\tilde\mu)=0$ and thus $\mu=\tilde\mu$.
$\Cox$

\bigskip
\textbf{Proof of Theorem \ref{Thm_Rate_Convergence}: }We show $c_N^{emp}(k,\nu')\to c(k,\nu')$ for all $k\in\{1,\dots,q\}$ and $\nu'\in\PP(\{1,\dots,q\})$.
For the nominator in the definition of $c_N^{emp}(k,\nu')$ we have
\begin{equation*}\label{Rat_Converence_Nominator}
\begin{split}
\mu_{\Phi,N-1}^{emp}[\nu']&(e^{-N\Phi(\frac{1}{N}\d_{\frac{2\pi}{q}
k}+\frac{N-1}{N}L_{N-1})+(N-1)\Phi(L_{N-1})})\cr
&=\frac{1}{Z_1}\int e^{-N\Phi(\frac{1}{N}\d_{\frac{2\pi}{q}
k}+\frac{N-1}{N}L_{N-1})}1_{T(L_{N-1})=\nu'}d\a^{\otimes N\ba i}\cr
&=\frac{1}{Z_1}\int e^{-N\Phi(L_{N-1})-d_{L_{N-1}}\Phi(\d_{\frac{2\pi}{q}
k}-L_{N-1})+o(\frac{1}{N})}1_{T(L_{N-1})=\nu'}d\a^{\otimes N\ba i}
\end{split}
\end{equation*}
where we used Taylor expansion. For the limit $N\uparrow\infty$ we can employ Varadhan's lemma together with Sanov's theorem and Lemma \ref{Fine_Descretization_General} and write
\begin{equation*}\label{Rat_Converence_Nominator_2}
\begin{split}
\frac{1}{Z_1}\int e^{-N\Phi(L_{N-1})-d_{L_{N-1}}\Phi(\d_{\frac{2\pi}{q}
k}-L_{N-1})+o(\frac{1}{N})}1_{T(L_{N-1})=\nu'}d\a^{\otimes N\ba i}
\to\frac{1}{Z_2}e^{-d_{\nu^{\nu'}}\Phi(\d_{\frac{2\pi}{q}
k}-\nu^{\nu'})}.\cr
\end{split}
\end{equation*}
Using the same arguments for the denominator of $c_N^{emp}(k,\nu')$ the constants cancel and we arrive at $c(k,\nu')$.
$\Cox$

\bigskip
\textbf{Proof of Theorem \ref{Thm_Rate_Convergence_Rotator}: }The first part of the theorem is an application of Lemma \ref{Fine_Descretization_General}. However we can use the special structur of the scalarproduct interaction to specify the constant $C(q(\Phi))$. Indeed, using the notation in the proof of  Lemma \ref{Fine_Descretization_General}, from \eqref{MF_EKO_key_1} we get
\begin{equation}\label{MF_EKO_key_2_Rotator}
\begin{split}
\Vert\g_k(\cdot|\mu)-\g_k(\cdot|\tilde\mu)\Vert\leq\frac{1}{4}\sup_{s,t\in S_k}|d_{\mu}\Phi(\d_s-\d_t)-d_{\tilde\mu}\Phi(\d_s-\d_t)|
\end{split}
\end{equation}
where $d_{\mu}\Phi(\d_s-\d_t)=-\b\langle \int\mu(d\o)e_\o,e_s-e_t\rangle$. With Cauchy-Schwartz 
\begin{equation}\label{MF_EKO_Estimate}
\begin{split}
\sup_{s,t\in S_k}&|\int(\tilde\mu(d\o)-\mu(d\o))\langle e_\o,e_s-e_t\rangle|\cr
&\leq\sup_{s,t\in S_k}\sup_{l\in\{1,\dots,q\}}|\int_{S_l}(\tilde\mu(d\o|S_l)-\mu(d\o|S_l))\langle e_\o,e_s-e_t\rangle|\cr
&\leq\sup_{l\in\{1,\dots,q\}}\sup_{s,t\in S_k}\sup_{x,y\in S_l}|\langle e_x-e_y,e_s-e_t\rangle|\Vert \tilde\mu(\cdot|S_l)-\mu(\cdot|S_l)\Vert\cr
&\leq4\sin^2\frac{\pi}{q}\sup_{l\in\{1,\dots,q\}}\Vert \tilde\mu(\cdot|S_l)-\mu(\cdot|S_l)\Vert.\cr
\end{split}
\end{equation}
By assumption $\b\sin^2(\frac{\pi}{q})<1$ and thus the first result follows.

Notice, in case of the standard scalarproduct potential we have $$-d_{\nu^{\nu'}}\Phi(\d_{\frac{2\pi}{q}
k}-\nu^{\nu'})=\b\langle\int\nu^{\nu'}(d\o)e_\o,e_{\frac{2\pi}{q}k}\rangle+\b\langle\int\nu^{\nu'}(d\o)e_\o,\int\nu^{\nu'}(d\o)e_\o\rangle$$
where the second summand is independent of the integration in the denominator of the rates and thus cancels. Using the notation $M_\b(\nu')=\int\nu^{\nu'}(d\o)e_\o$ we arrive at definition of the rates \eqref{Rotation_Flow_Rates}.
$\Cox$

\bigskip
To complement the above criterion on the finess of descretization in order to have unique constrained free energy minimizers for the rotator model, let us consider an equivalent of a checkerboard configuration on the lattice. Namely the measure with equal weight on segments facing in opposite directions. This will lead to a criterion for non-uniqueness of the constrained minimizers. For convenience take $q$ even. We condition on $\nu'=\frac{1}{2}(\d_1+\d_{\frac{q}{2}+1})$, then from \eqref{MF_DLR_General} we know for a constrained minimizers $\nu^{\nu'}$ we have
\begin{equation*}\label{Complement_Criterion}
\begin{split}
M_\b(\nu')=\sum_k\nu'(k)\int_{S_k}\nu^{\nu'}(d\o|S_k)e_\o=\sum_k\nu'(k)\frac{\int_{S_k}e_\o\exp(\b\langle e_{\o},M_\b(\nu')\rangle)\a(d\o)}{\int_{S_k}\exp(\b\langle e_{\o},M_\b(\nu')\rangle)\a(d\o)}.
\end{split}
\end{equation*}
Note, this equation is often referred to the mean-field equation. By symmetry and under suitable coordinates this fixed point equation becomes one-dimensional and reads
\begin{equation}\label{Non_Gibbs_MFE}
\begin{split}
m=
\frac{\int_{-\frac{\pi}{q}}^{\frac{\pi}{q}}\sin(\o)\exp(\b m \sin(\o))\a(d\o)}{\int_{-\frac{\pi}{q}}^{\frac{\pi}{q}}\exp(\b m \sin(\o))\a(d\o)}=:F_{q}(\b m).
\end{split}
\end{equation}
%
Since $F_{q}$ is concave, the equation \eqref{Non_Gibbs_MFE} has no non-trivial fixed point if $\frac{d}{d m}_{|m=0}F_{q}(m)<1/\b$, i.e. $\b(\frac{1}{2}-\frac{q}{4\pi}\sin(\frac{2\pi}{q}))<1$. On the other hand if
\begin{equation}\label{Non_Gibbs_Criterion}
\begin{split}
\b(1-\frac{q}{2\pi}\sin(\frac{2\pi}{q}))>2
\end{split}
\end{equation}
 there must be a non-trivial fixed point since $F$ is bounded and continuous. In other words if \eqref{Non_Gibbs_Criterion} holds, there are two distinct measures $\nu^{\nu'}_+\neq\nu^{\nu'}_-$. In particular $\int\nu^{\nu'}_+(d\o)e_\o=-\int\nu^{\nu'}_-(d\o)e_\o\neq0$ and $\Psi(\nu^{\nu'}_+)=\Psi(\nu^{\nu'}_-)$ because of symmetry. Hence in the regime \eqref{Non_Gibbs_Criterion} we just provided an example were the constrained model has multiple Gibbs measures. 

\subsection{Infinite-volume rotation dynamics}\label{The infinite-volume rotation dynamics}

Let us in the sequel specify to the rotator model with scalarproduct potential and its discretization, assumed to be in the parameter regime $\b>2$ and $\b\sin^2(\frac{\pi}{q})<1$. 
\begin{lem}\label{ODE_Solvable}
The non-linear system of ordinary differential equations given in Theorem \ref{Path_LLN} with rates \eqref{Rotation_Flow_Rates}
is uniquely solvable globally in time.
\end{lem}
Notice \eqref{Rotation_Flow_RHS} can be interpreted as inflow from below into state $k$ minus outflow in the direction $k+1$. 

\bigskip
\textbf{Proof of Lemma \ref{ODE_Solvable}: }For a given initial measure $\nu'_0$ the system \eqref{Rotation_Flow} is uniquely solvable locally in time by the Picard-Lindel\" of theorem. Indeed, we show Lipschitz continuity of \eqref{Rotation_Flow_RHS} as a function of $\nu'$ w.r.t the total-variation distance. It suffices to show Lipschitz continuity for $\nu'\mapsto \nu^{\nu'}$ since \eqref{Rotation_Flow_RHS} is a composition of Lipschitz continuous functions of $\nu^{\nu'}$.
First note
\begin{equation*}\label{MF_EKO_6}
\begin{split}
\Vert \nu^{\nu'}-\nu^{\tilde\nu'}\Vert&\leq\sup_{k\in\{1,\dots,q\}}\Vert\nu^{\nu'}(\cdot|S_k)-\nu^{\tilde\nu'}(\cdot|S_k)\Vert+\Vert\nu'-\tilde\nu'\Vert.\cr
\end{split}
\end{equation*}
Introducing $\g_k(ds|\nu^{\nu'})$ as defined in \eqref{MF_DLR_General} we can further write for a bounded measurable function $f$
\begin{equation*}\label{Coarse_Grained_Specification}
\begin{split}
|\nu^{\nu'}(f|S_k)-\nu^{\tilde\nu'}(f|S_k)|&=|\g_k(f|\sum_{l=1}^q\nu'(l)\nu^{\nu'}(\cdot|S_l))-\g_k(f|\sum_{l=1}^q\tilde\nu'(l)\nu^{\tilde\nu'}(\cdot|S_l))|\cr
&\leq \d(f)\b\sin^2(\frac{\pi}{q})\sup_{l\in\{1,\dots,q\}}\Vert\nu^{\nu'}(\cdot|S_l)-\nu^{\tilde\nu'}(\cdot|S_l)\Vert\cr
&\hspace{1cm}+\Vert\g_k(\cdot|\sum_{l=1}^q\nu'(l)\nu^{\tilde\nu'}(\cdot|S_l))-\g_k(\cdot|\sum_{l=1}^q\tilde\nu'(l)\nu^{\tilde\nu'}(\cdot|S_l))\Vert)\cr
\end{split}
\end{equation*}
where we used \eqref{MF_EKO_key_2_Rotator} and \eqref{MF_EKO_Estimate} for the first summand. For the second summand we have
\begin{equation*}\label{MF_EKO_Estimate_2}
\begin{split}
\Vert\g_k&(\cdot|\sum_{l\in\{1,\dots,q\}}\nu'(l)\nu^{\tilde\nu'}(\cdot|S_l))-\g_k(\cdot|\sum_{l\in\{1,\dots,q\}}\tilde\nu'(l)\nu^{\tilde\nu'}(\cdot|S_l))\Vert)\cr
&\leq \frac{\b}{4}\sup_{s,t\in S_k}|\langle \sum_{l=1}^q\nu'(l)\int\nu^{\tilde\nu'}(\cdot|S_l)(d\o)e_\o-\sum_{l=1}^q\tilde\nu'(l)\int\nu^{\tilde\nu'}(\cdot|S_l)(d\o)e_\o,e_s-e_t\rangle|\cr
&\leq \frac{\b}{4}\sup_{s,t\in S_k}\sup_{i,j\in\{1,\dots,q\}}|\langle \int\nu^{\tilde\nu'}(\cdot|S_i)(d\o)e_\o,-\int\nu^{\tilde\nu'}(\cdot|S_j)(d\o)e_\o,e_s-e_t\rangle|\Vert\nu'-\tilde\nu'\Vert\cr
&\leq \frac{\b}{2}\sup_{s,t\in S_k}\Vert e_s-e_t\Vert\Vert\nu'-\tilde\nu'\Vert=\b\sin(\frac{\pi}{q})\Vert\nu'-\tilde\nu'\Vert\cr
\end{split}
\end{equation*}
where we used \eqref{MF_EKO_key_1} in the first inequality. Thus taking the supremum over $f$ and $k$ and using the fact, that we are in the right parameter regime, we have 
\begin{equation*}\label{MF_EKO_Goergiis_Argument}
\begin{split}
&\sup_{k\in\{1,\dots,q\}}\Vert\nu^{\nu'}(\cdot|S_k)-\nu^{\tilde\nu'}(\cdot|S_k)\Vert\leq C \Vert\nu'-\tilde\nu'\Vert.
\end{split}
\end{equation*}
But this is Lipschitz continuity.

Solutions also always exist globally: If $\nu'_t(k)=0$ for some $k$, we have $\frac{d}{dt}\nu'_t(k)=c(k-1,\nu'_t)\nu'_t(k-1)\geq0$. In other words, if a solution is on the boundary of the simplex, the vector field forces the trajectory back inside the simplex.
$\Cox$

\bigskip
\textbf{Remark: }The above lemma in particular proves, that the so called \textit{second-layer mean-field specification} $\g'(k|\nu'):=\frac{\int_{S_k}\exp(\b\langle M_\b(\nu'),e_\o\rangle)\a(d\o)}{\int\exp(\b\langle M_\b(\nu'),e_\o\rangle)\a(d\o)}$ is continuous w.r.t the boundary entry $\nu'$. This is the defining property for a system after coarse-graining to be called Gibbs.
%



\bigskip
\textbf{Proof of Theorem \ref{Path_LLN}: }
We use the Feng-Kurtz scheme as presented in \cite{FK06,EFHR10,RW12} to show convergence on the level of trajectories. The Feng-Kurtz Hamiltonian for the generator $Q_N^{emp}$ reads
\begin{equation*}\label{rotation_generator4}
\begin{split}
\mathcal{H}(\nu',f)&=\sum_{k=1}^q\nu'(k)c(k,\nu')(e^{df_{\nu'}(\d_{k+1}-\nu')-df_{\nu'}(\d_{k}-\nu')}-1)
\end{split}
\end{equation*}
where $f: \PP(\{1,\dots,q\})\to\R$ is a differentiable observable and we used the convergence of the rates from Theorem \ref{Thm_Rate_Convergence_Rotator}.
This Hamiltonian is of the form as presented in \cite{FK06} Section 10.3. with $b(\nu'):=F(\nu')$ and $\eta(\nu',\d_{k+1}-\d_k):=\nu'(k)c(k,\nu')$. Following the roadmap of \cite{FK06} we verify (using references as in \cite{FK06}):
\begin{enumerate}
\item $X_n$ is exponentially tight in the path space by Theorem 4.1. since
\begin{equation*}\label{FK_Theorem41}
\begin{split}
\E(e^{N\l\Vert X_{t+\d}^N-X_{t}^N\Vert}|\mathcal{F}_t^N)\leq\E(e^{\l \text{Pois}(N K\d)})=e^{N K\d(e^{\l}-1)}
\end{split}
\end{equation*}
where $\text{Pois}(r)$ stands for a Poisson random variable with intensity $r$, $K:=\sup_{\nu'\in\\P(\{1,\dots,q\}),N\in\N,k\in\{1,\dots,q\}}c^{emp}_N(k,\nu')$ and $\mathcal{F}_t^N=\s((X_s)^N_{0\leq s\leq t})$. In particular $\lim_{\d\downarrow0}K\d(e^{\l}-1)=0$ and thus criterion b) of Theorem 4.1. is satisfied.
\item The comparison principle holds for the generator $\mathcal{H}$ (ensuring the existence of the so-called exponential semigroup corresponding to $\mathcal{H}$) since 
the conditions of Lemma 10.12. are satisfied. In particular we used the Lipschitz continuity of $F(\nu')$ from our Lemma \ref{ODE_Solvable}. This property ensures existence of a LDP for the finite-dimensional distributions of the process.
\end{enumerate}
Further since also (10.18) in \cite{FK06} is satisfied, Theorem 10.17. in \cite{FK06} gives the LDP on the level of paths implicitly via the exponential semigroup. 

In order to get a nice variational representation of the large deviation rate function we calculate the Lagrangian of $\mathcal{H}$ for the measure $\nu'$ and the velocity $u'$ (a zero-weight signed measure on $\{1,\dots,q\}$) as in Lemma 10.19. in \cite{FK06}
\begin{equation*}\label{Lagrangian_0}
\begin{split}
L(\nu',u')&:=\sup_p\bigl(\langle p,u'\rangle -\sum_{k=1}^q\nu'(k)c_L(k,\nu')(e^{p(k+1)-p(k)}-1)\bigr)\geq 0.\cr
\end{split}
\end{equation*}
Let $\mathbb{P}^{N}_{{\nu_0'}}$ denote the law of the Markov process $Q_N^{emp}$ started in $\nu'_0$, then by Theorem 10.22. in \cite{FK06} we have
\begin{equation*}\label{Path_LDP}
\begin{split}
\mathbb{P}^{N}_{{\nu'_0}}\Bigl((\cdot_t)_{t\in[0,T_f]}\approx(\nu'_t)_{t\in[0,T_f]}\Bigr)\approx\exp(-N\int_{0}^{T_f}L(\nu'_t,\frac{d}{dt}\nu'_t)dt)
\end{split}
\end{equation*}
where 
the approximation signs should be understood in the sense of the LDP with the Skorokhod topology on the space of cadlag paths. In fact by \cite{FK06} Theorem 4.14 the LDP even holds in the uniform topology. 

To obtain the LLN we need to show $L(\nu'_t,\frac{d}{dt}\nu'_t)=0$ if $\frac{d}{dt}\nu'_t$ is given by \eqref{Rotation_Flow}. But this is true: The Lagrangian for \eqref{Rotation_Flow} reads
\begin{equation}\label{Lagrangian}
\begin{split}
L&(\nu'_t,\frac{d}{dt}\nu'_t)=\sup_p\bigl(\langle p,\frac{d}{dt}\nu'_t\rangle -\sum_{k=1}^q\nu'_t(k)c(k,\nu'_t)(e^{p(k+1)-p(k)}-1)\bigr)\cr
&=\sup_p\bigl(\sum_{k=1}^q\nu'_t(k)c(k,\nu'_t)(p(k+1)-p(k)+1-e^{p(k+1)-p(k)})\bigr)=:\sup_pJ(p)\cr
\end{split}
\end{equation}
with $\frac{\partial}{\partial p(l)}J(p)=\nu'_t(l-1)c(l-1,\nu'_t)(1-e^{p(l)-p(l-1)})-\nu'_t(l)c(l,\nu'_t)(1-e^{p(l+1)-p(l)})$.
In case $\nu'_t(l)=0$ for some $l\in\{1,\dots,q\}$ and $p^*$ realizing the supremum in \eqref{Lagrangian} we have $\nu'_t(k)c(k,\nu'_t)(1-e^{p^*(k+1)-p^*(k)})=0$ for all $k\in\{1,\dots,q\}$ and since $c(k,\nu'_t)>0$ in particular $p^*(k+1)=p^*(k)$ whenever $\nu'_t(k)>0$. Thus $L(\nu'_t,\frac{d}{dt}\nu'_t)=0$. In case $\nu'_t(k)>0$ for all $k\in \{1,\dots,q\}$, $J$ is strictly concave away from any constant vector $p=(c,\dots,c)^T\in\R^q$, to be precise
\begin{equation*}\label{Lagrangian_6}
\begin{split}
\langle z,(\frac{\partial^2}{\partial p(i)\partial p(j)}J)_{i,j} z\rangle=-\sum_{k=1}^q\nu'_t(k)c(k,\nu'_t)e^{p(k+1)-p(k)}(z_{k+1}-z_k)^2
\end{split}
\end{equation*}
for all $q$-dimensional vectors $z$, and thus $J(0)$ is the global maximum 
of $J$. Hence $L(\nu'_t,\frac{d}{dt}\nu'_t)=0$. Since $L(\nu'_t,\cdot)$ is strictly convex as a Legendre tranform of the strictly convex Feng-Kurtz Hamiltonian $\mathcal{H}$, the flow \eqref{Rotation_Flow} is the unique dynamics such that $
\int_0^{T_f}L(\nu'_t,\frac{d}{dt}\nu'_t)dt=0$. But that means, according to the LDP, that \eqref{Rotation_Flow} is the unique limiting dynamics as $N\to\infty$.
$\Cox$

\bigskip
\textbf{Remark: }
If one is only interested in the (weak) LLN for $X_n$, one can also apply Theorem 2 of \cite{O84}, with the minor alteration, that our rates are $N$-dependent but convergent. The proof of the result uses martingale respesentation to derive tightness of $(\mathbb{P}^{N}_{{\nu_0'}})_{N\in\N}$
(in the pathspace equipped with the Skorokhod topology.) The uniqueness of the limiting (deterministic) process is shown by a coupling argument. 
For the sake of accessibility we compare the notation in \cite{O84} with ours: $\g_\o(k):=\d_1(k)$, $A(k,y,\nu'):=c_N^{emp}(k,\nu')$, all given topologies on $\PP(\{1,\dots,q\})$ 
are equivalent, the Lipschitz condition (B4) is satisfied since $\nu'\mapsto c_N^{emp}(k,\nu')$ is a composition of Lipschitz continuous functions (where we have to use Lemma \ref{ODE_Solvable}), (B3), (B2) and (B1) are trivially satisfied.

%

\section{Properties of the flow}\label{Properties of the flow}
In \cite{JK12} one of the main results states the existence of a unique translation-invariant invariant measure for the rotation dynamics combined with a Glauber dynamics on the lattice, which is not long-time limit of all starting measures. This is done by identifying a set of starting measures, namely the set of extremal translation-invariant Gibbs measures of a discretized version of the XY model, which is not attracted to the invariant measure. Further results about attractivity of the rotation dynamics alone for general starting measures seemed to be difficult on the lattice. 

In this section we reproduce the equivariance properties of the discretization map for the dynamical system. Further we investigate attractivity properties of the flow.

Before we start let us note, that (as in the lattice situation) the commutator (in the form of the Lie bracket) of the rotation dynamics and the corresponding Glauber dynamics vanishes on $\GG'$. In general there is no reason to believe that the two dynamics do commute.

In the sequel we denote $\nu_t'$ a probability measure on $\{1,\dots,q\}$ at time $t$ under the rotation dynamics \eqref{Rotation_Flow}. 
\subsection{Closed orbit and equivariance of the discretization map}\label{Properties of the flow: Closed orbit and equivariance property of the discretization map}

\textbf{Proof of Theorem \ref{Diagram_Commutating}: }For $\nu\in\GG(\Phi)$ we have $T(\nu)\in\GG'$ by the contraction principle. Further for $\nu'\in\GG'$ we have
\begin{equation*}\label{Gibbs_back}
\begin{split}
\inf_{\nu}\Psi(\nu)=\inf_{\tilde\nu'}\inf_{\nu:T\nu=\tilde\nu'}\Psi(\nu)
=\inf_{\tilde\nu'}\Psi'(\tilde\nu')=\Psi'(\nu')=\Psi(\nu^{\nu'}),
\end{split}
\end{equation*}
hence $\nu^{\nu'}\in\GG(\Phi)$ and we have established a one-to-one correspondence between $\GG'$ and $\GG({\Phi})$.

Let us verify the dynamical aspects of the diagram.
Let $\nu_t\in\GG(\Phi)$ 
and compute the derivative (in analogy to Proposition \ref{Prop_Finite_Generator_Time_Dependent}) and note that indeed the left-sided and the right-sided derivatives coincide
\begin{equation}\label{Infinitesimal_Rotation}
 \begin{split}
\frac{d}{d\e}|_{\e=0}&T(\nu_{t+\e})(k)=\frac{d}{d\e}|_{\e=0}\frac{\int_{S_k}\a(d\s)\exp(\b\langle e_\s,M_\b(T(\nu_{t+\e}))\rangle)}{
\int\a(d\s)\exp(\b\langle e_\s,M_\b(T(\nu_{t+\e}))\rangle)}\cr
&=\frac{d}{d\e}|_{\e=0}\frac{\int_{\frac{2\pi}{q}
(k-1)-\e}^{\frac{2\pi}{q}
k-\e}\a(d\s)\exp(\b\langle e_\s,M_\b(T(\nu_{t}))\rangle)}{
\int\a(d\s)\exp(\b\langle e_\s,M_\b(T(\nu_{t}))\rangle)}\cr
&=\frac{\exp(\b\langle e_{\frac{2\pi}{q}
(k-1)},M_\b(T(\nu_{t}))\rangle)-\exp(\b\langle e_{\frac{2\pi}{q}
k},M_\b(T(\nu_{t}))\rangle)}{
\int\a(d\s)\exp(\b\langle e_\s,M_\b(T(\nu_{t}))\rangle)}\cr
&=c(k-1,T(\nu_t))T(\nu_t)(k-1)-c(k,T(\nu_t))T(\nu_t)(k)=F(T(\nu_t))(k).
\end{split}
\end{equation}

By Lemma \ref{ODE_Solvable}, the differential equation \eqref{Rotation_Flow} is uniquely solvable globally in time. 
Since $(T(\nu_t))_{t\geq 0}$ is a trajectory in $\PP(\{1,\dots,q\})$ solving the differential equation we have $T(\nu_{t+s})=\phi(s,T(\nu_t))$ for all $s,t\geq0$.

\bigskip
Note: One can also 
show higher differentiability of the flow with respect to the initial condition. Strong enough differentiability of $F$ would ensure that. This again would be guaranteed by strong enough differentiability of $\nu'\mapsto M_\b(\nu')$. One can employ an implicit function theorem applied to the mean-field equation to get that kind of regularity. Unfortunately a 
price to pay could a priori be the assumption of an unspecified maybe large $q$, so some additional 
technical work would be needed. 
$\Cox$

\bigskip
In the sequel we will often refer to $\GG'$ as the\textit{ periodic orbit} of the flow $(\phi_t)_{t\geq0}$. 

\subsection{Attractivity of the closed orbit via free energy}\label{Attractivity via the free energy}

\begin{lem}\label{Time_Derivative_Lyapunov}
The time derivative of the free energy on $\PP(\{1,\dots,q\})$ reads
\begin{equation}\label{Lyapunov_Derivative}
\begin{split}
&\frac{d}{dt}|_{t=0}\Psi'(\phi(t,\nu'))\cr
&=\sum_{k\in\{1,\dots,q\}} \frac{e^{\b\langle e_{\frac{2\pi}{q}
k},M_\b(\nu')\rangle}}{\int_{S_k} e^{\b\langle e_\o,M_\b(\nu')\rangle}\a(d\o)}
\nu'(k)\log\frac{\nu'(k+1)\int_{S_k} e^{\b\langle e_\o,M_\b(\nu')\rangle}\a(d\o)}{\nu'(k)\int_{S_{k+1}} e^{\b\langle e_\o,M_\b(\nu')\rangle}\a(d\o)}.\cr
\end{split}
\end{equation}
\end{lem}
\textbf{Proof of Lemma \ref{Time_Derivative_Lyapunov}: }We have $\Psi'(\nu')=\Psi(\nu^{\nu'})=S(\nu^{\nu'}|\a)+\Phi(\nu^{\nu'})+Const$
where $\Phi(\nu^{\nu'})=-\frac{\b}{2}\Vert \int\nu^{\nu'}(d\o)e_{\o}\Vert^2=-\frac{\b}{2}\Vert M_\b({\nu'})\Vert^2$ and
\begin{equation*}\label{Relative_Entropy_Constrained_Mini}
\begin{split}
&S(\nu^{\nu'}|\a)=\int d\nu^{\nu'}\log\frac{d\nu^{\nu'}}{d\a}=\sum_{k\in\{1,\dots,q\}}\nu'(k)\int\nu^{\nu'}(d\o|S_k)\log\frac{\nu'(k)e^{\b\langle e_\o,M_\b(\nu')\rangle}}{\int_{S_k} e^{\b\langle e_{\tilde\o},M_\b(\nu')\rangle}\a(d\tilde\o)}\cr
&=\sum_{k\in\{1,\dots,q\}}\nu'(k)\int\nu^{\nu'}(d\o|S_k)[\log\nu'(k)+\b\langle e_\o,M_\b(\nu')\rangle-\log\int_{S_k}e^{\b\langle e_\o,M_\b(\nu')\rangle}\a(d\o)]\cr
&=S(\nu'|\a')+\b\Vert M_\b(\nu')\Vert^2-
\sum_k\nu'(k)\log (q\int_{S_k}e^{\b\langle e_\o,M_\b(\nu')\rangle}\a(d\o))\cr
\end{split}
\end{equation*}
where we used $\nu^{\nu'}=\sum_{k\in\{1,\dots,q\}}\nu'(k)\int\nu^{\nu'}(\cdot|S_k)$ and wrote $\a':=\frac{1}{q}\sum_{k=1}^q\d_k$ for the equidistribution. The time derivative is now just a simple calculation.
$\Cox$

\bigskip
\textbf{Proof of Proposition \ref{Lyapunov_Lemma}: }First note, if $\nu'\in\GG'$ or $\nu'=\frac{1}{q}\sum_{k=1}^q\d_k$
we have $\nu'(k)=K\int_{S_k}e^{\b\langle e_\o,M_\b(\nu')\rangle}\a(d\o)$ and hence $\frac{d}{dt}|_{t=0}\Psi'(\phi(t,\nu'))=0$.
For any distribution with no weight on at least one $k\in \{1,\dots,q\}$, the r.h.s of \eqref{Lyapunov_Derivative} is minus infinity. 

Let us change the perspective and assume $\b M_\b(\nu')=:x\in\R^2$ to be given instead of $\nu'$. Let $\G_k: \R^2\to\R^2, \G_k(x_1,x_2)=\frac{\int_{S_k}e_\o\exp(\langle e_\o,(x_1,x_2)\rangle)\a(d\o)}{\int_{S_k}\exp(\langle e_\o,(x_1,x_2)\rangle)\a(d\o)}$ and define $\L_x:=\{\{\l_1, ...,\l_q\}\in\R^q_+|\sum_k\l_k\G_k(x)=x\}$ the space of unnormalized measures such that their corresponding probability measures $(\l_k/\Vert \l\Vert_{l^1})_{k\in\{1,...,q\}}$ have magnetization $x/\b$,
in particular $\Vert \l\Vert_{l^1}=\b$.
Let us rewrite the free energy and prove instead of $\frac{d}{dt}|_{t=0}\Psi'(\phi(t,\nu'))\leq 0$,
\begin{equation*}\label{Second_Layer_Free_Energy_Time_Derivative_New_Perspective}
\begin{split}
0&\geq\sum_{k\in\{1,\dots,q\}} \frac{e^{\langle e_{\frac{2\pi}{q}
k},x\rangle}}{\int_{S_k} e^{\langle e_\o,x\rangle}\a(d\o)}
\l(k)\log\frac{\l(k+1)\int_{S_k} e^{\langle e_\o,x\rangle}\a(d\o)}{\l(k)\int_{S_{k+1}} e^{\langle e_\o,x\rangle}\a(d\o)}=:G_x(\l_1, ...\l_q)\cr
\end{split}
\end{equation*}
for $\l\in\L_x$. One way to do this is to show that for given $x\in\R^2$ the maximum of $G_x$ under the constraint $\{\l_1, ...\l_q\}\in\L_x$ is lower or equal zero. Let us apply Lagrange multipliers $\a_1,\a_2$, then we must solve the following $q+2$ equations
\begin{equation}\label{Lagrange_equations}
\begin{split}
\frac{d}{d\l_k}[G_x(\l_1, ...\l_q)+\a_1(\sum_k\l_k\G_{k}(x)_1-x_1)+\a_2(\sum_k\l_k\G_{k}(x)_2-x_2)]&=0\cr
\sum_k\l_k\G_{k}(x)&=x.
\end{split}
\end{equation}
The first line of \eqref{Lagrange_equations} reads 
\begin{equation*}\label{Second_Layer_Free_Energy_5}
\begin{split}
\frac{e^{\langle e_{\frac{2\pi}{q}k},x\rangle}}{\int_{S_k} e^{\langle e_\o,x\rangle}\a(d\o)}\log\frac{\l(k+1)\int_{S_k} e^{\langle e_\o,x\rangle}\a(d\o)}{\l(k)\int_{S_{k+1}} e^{\langle e_\o,x\rangle}\a(d\o)}-\frac{e^{\langle e_{\frac{2\pi}{q}k},x\rangle}}{\int_{S_k} e^{\langle e_\o,x\rangle}\a(d\o)}&\cr
+\frac{e^{\langle e_{\frac{2\pi}{q}(k-1)},x\rangle}}{\int_{S_{k-1}} e^{\langle e_\o,x\rangle}\a(d\o)}
\frac{\l(k-1)}{\l(k)}+\a_1\G_{k}(x)_1+\a_2\G_{k}(x)_2&=0.
\end{split}
\end{equation*}
Multiplying these equations with $\l_k$, summing and applying the constraint condition we have $G_x(\l_1, ...\l_q)+\langle \a,x\rangle=0$
and thus the sign is determined by the Lagrange multipliers. 
Define $\b>0$ such that $\frac{x}{\b}=\frac{\int e_\o e^{\langle e_\o,x\rangle}\a(d\o)}{\int e^{\langle e_\o,x\rangle}\a(d\o)}$ and $\l_k:=\b\frac{\int_{S_k} e^{\langle e_\o,x\rangle}\a(d\o)}{\int e^{\langle e_\o,x\rangle}\a(d\o)}$ the Gibbs measure for $x$ rescaled with $\b$. In particular $\sum_k\l_k\G_k(x)=x$. First we show $\underline{\l}:=(\l_i)_{i\in\{1,\dots,q\}}$ is an extremal point of $G_x$ under the constraint $\L_x$. Indeed, set $\a_1=x_2$ and $\a_2=-x_1$, then the first $q$ equations in \eqref{Lagrange_equations} read
\begin{equation*}\label{Gibbs_Extremal}
\begin{split}
\frac{e^{\langle e_{\frac{2\pi}{q}(k-1)},x\rangle}-e^{\langle e_{\frac{2\pi}{q}k},x\rangle}}{\int_{S_k} e^{\langle e_\o,x\rangle}\a(d\o)}+x_2\frac{\int_{S_k} \cos(\o) e^{\langle e_\o,x\rangle}\a(d\o)}{\int_{S_k} e^{\langle e_\o,x\rangle}\a(d\o)}-x_1\frac{\int_{S_k} \sin(\o) e^{\langle e_\o,x\rangle}\a(d\o)}{\int_{S_k} e^{\langle e_\o,x\rangle}\a(d\o)}.
\end{split}
\end{equation*}
But this is zero since $\int_{S_k} (x_2\cos(\o)-x_1\sin(\o))e^{\langle e_\o,x\rangle}\a(d\o)=e^{\langle e_{\frac{2\pi}{q}(k-1)},x\rangle}-e^{\langle e_{\frac{2\pi}{q}k},x\rangle}$.

\bigskip
Secondly we show $G_x$ is concave on $\L_x$, indeed 
\begin{equation*}\label{G_Derivative}
\begin{split}
\frac{\partial}{\partial \l_k}G_x(\underline\l)=c(k,x)\log\frac{\l(k+1)\int_{S_k} e^{\langle e_\o,x\rangle}\a(d\o)}{\l(k)\int_{S_{k+1}} e^{\langle e_\o,x\rangle}\a(d\o)}-c(k,x)
+c(k-1,x)\frac{\l(k-1)}{\l(k)}
\end{split}
\end{equation*}
thus the Hessian matrix has non-zero entries only on the diagonal and on the two neighboring diagonals 
\begin{equation*}\label{G_Hessian}
\begin{split}
(\frac{\partial}{\partial \l_k})^2G_x(\underline\l)&=-(\frac{c(k,x)}{\l_k}+\frac{c(k-1,x)\l_{k-1}}{\l_k^2})\cr
\frac{\partial^2}{\partial\l_{k+1}\partial \l_k}G_x(\underline\l)&=\frac{c(k,x)}{\l_{k+1}}\hspace{0.5cm}\text{ and }\hspace{0.5cm}\frac{\partial^2}{\partial\l_{k-1}\partial \l_k}G_x(\underline\l)=\frac{c(k-1,x)}{\l_{k}}.\cr
\end{split}
\end{equation*}
In order to check definiteness we apply an arbitrary vector $(\l_1z_1, \dots,\l_qz_q)^T$ from both sides, which gives us  
\begin{equation*}\label{Hessian_Definiteness}
\begin{split}
\sum_{i=1}^q(w_iz_iz_{i+1}+w_{i-1}z_{i-1}z_{i}-(w_i+w_{i-1})z_i^2)=\sum_{i=1}^qw_i(2z_iz_{i+1}-(z_i^2+z_{i+1}^2))
\end{split}
\end{equation*}
where we wrote $w_i:=\l_ic(i,x)$. Since $z_iz_{i+1}\leq\frac{z_i^2+z_{i+1}^2}{2}$ the Hessian is negative semidefinite and thus $G_x$ is concave. Hence $\underline\l$ must be a global maximum for $G_x$. Notice, the eigenspace for the eigenvalue zero is $\{v\in\R^q|v=c\underline\l\}$. Thus the only direction in which $D^2G_x(\underline\l)$ is not strictly negative is the one along $\underline\l$, but $c\underline\l\notin\L_x$ unless $c=1$. Hence $\underline\l$ is the only maximum in $G_x$. Since all $G_x$ are disjoint and every probability measure belongs to some $G_{x/\b}$, we showed that $\Psi'$ is indeed strictly decreasing away from the peridic orbit and the equidistribution.
$\Cox$

\begin{figure}[h]
\begin{center}
\includegraphics[width=7.2cm]{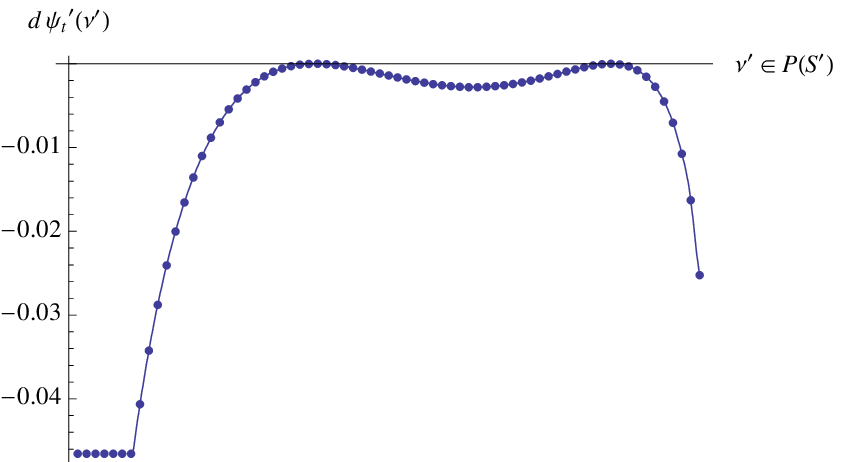}
\includegraphics[width=7.2cm]{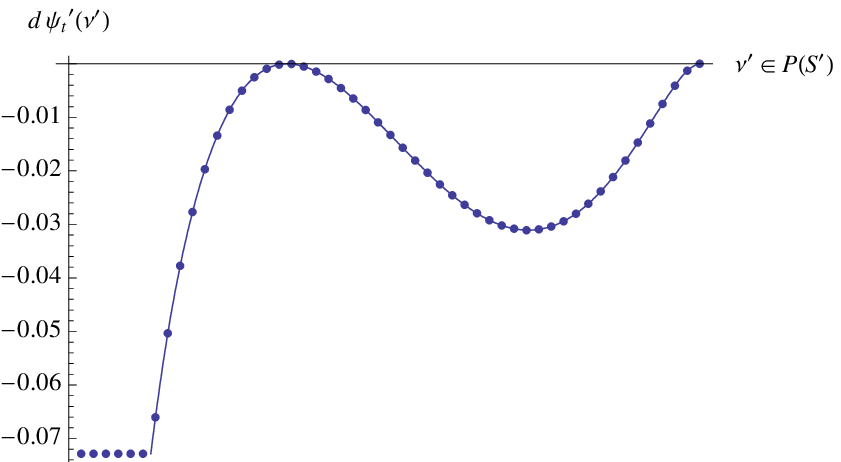}
\end{center}
\caption{\scriptsize{Numerical simulations for the time-derivative of the Lyapunov function along the line (of ten-dimensional probability measures) going through the equidistribution and the Gibbs measure with order parameter $m^*e_0$. The left figure (respectively on the right) shows the situation with $\b=2.3$ (respectively $\b=3$.) In both figures the left local maximum belongs to the equidistribution. The right local maximum is at the Gibbs measure. In all cases, at the boundaries of the simplex, the function goes to minus infinity. $S'=\{1,\dots,q\}$.}}
\end{figure}

%
%

\textbf{Proof of Theorem \ref{Attract_Result}: } 
Let 
$K_\e:=\{\nu'|\inf_{\mu'\in\GG'}\sum_{k=1}^q|\nu'(k)-\mu'(k)|\geq\e\}$ then $\frac{d}{ds}|_{s=0}\Psi'(\phi(s,\nu'))\leq\d_\e<0$ for all $\nu'\in K_\e$ for some $\d_\e<0$ by compactness of $K_\e$ and Proposition \eqref{Lyapunov_Lemma}. Assume $(\phi(t,\nu'))_{t\geq0}\in K_\e$ for all $t\geq0$, then for all $s\geq0$
\begin{equation*}\label{Compactness_Argument}
\begin{split}
\Psi'(\phi(t+s,\nu'))-\Psi'(\phi(t,\nu'))=\int_0^s\frac{d}{dh}\Psi'(\phi(t+h,\nu'))dh\leq\d_\e s.
\end{split}
\end{equation*}
But
$\Psi'\geq0$, 
a contradiction. Since $\Psi'(\nu')<\Psi'(\frac{1}{q}\sum_k\d_k)$ for any Gibbs measure $\nu'$ and by assumption $\Psi'(\nu')<\Psi'(\frac{1}{q}\sum_k\d_k)$, $\phi(t,\nu')$ can only leave $K_\e$ towards the periodic orbit.
$\Cox$

\subsection{Local stability analysis at the equidistribution via linearization}\label{Local stability analysis of the equidistribution via linearization}
Recall the definition of the flow
\begin{equation*}\label{Rotation_Flow_2}
\begin{split}
\frac{d}{dt}|_{t=0}\phi(t,\nu')(k)=c(k-1,\nu')\nu'(k-1)-c(k,\nu')\nu'(k)=F(\nu')(k).
\end{split}
\end{equation*}
In order to understand local attractivity, we calculate the linearized r.h.s $dF$. To simplify notation, let us write just $\frac{d}{d\e}$ when we mean $\frac{d}{d\e}|_{\e=0}$, $\nu^{\nu'}_k:=\nu^{\nu'}(\cdot|S_k)$ and $m_\b(\nu):=\int\nu(d\o)e_\o$. For any $\nu'\in \PP(\{1,\dots,q\})$ and zero-weight signed measure $\r$ on $\{1,\dots,q\}$ we have
\begin{equation*}\label{Linearization_RHS}
 \begin{split}
d F_{\nu'}(\r)(k)&=\frac{d}{d\eps}F(\nu' + \eps \rho)(k)=
( \frac{d}{d\eps}c(k-1,\nu' + \eps \rho))\nu'(k-1)\cr
&\hspace{0,5cm}- (\frac{d}{d\eps}c(k,\nu' + \eps \rho))\nu'(k)+c(k-1,\nu')\rho(k-1)- c(k,\nu')\rho(k)\cr
\end{split}
\end{equation*}
with
\begin{equation*}
 \begin{split}
&\frac{d}{d\e}c(k,\nu' + \eps \rho)=\b c(k,\nu')\langle e_{\frac{2\pi}{q}k}
-\int_{S_k} e^{\b\langle e_\o,M_\b(\nu')\rangle}
e_\o\a(d\o),\frac{d}{d\e}M_\b(\nu'+\e\r)\rangle\cr
&=\b [c(k,\nu')\langle e_{\frac{2\pi}{q}k}
,\frac{d}{d\e}M_\b(\nu'+\e\r)\rangle-e^{\b\langle e_{\frac{2\pi}{q}k},M_\b(\nu')\rangle}
\langle m_\b(\nu_k^{\nu'}),\frac{d}{d\e}M_\b(\nu'+\e\r)\rangle].\cr
\end{split}
\end{equation*}
To compute the derivative $\frac{d}{d\eps}M_\b(\nu' + \eps \rho)$, we use the $2q$-dimensional mean-field equation and apply the implicit function theorem. We have
\begin{equation*}\label{Local_Magnetisation_1}
 \begin{split}
\frac{d}{d\e}M_\b(\nu'+\e\r)&=\sum_k \nu'(k)\frac{d}{d\e}m_\b(\nu^{\nu'+\e\r}_k)+\sum_k\r(k)m_\b(\nu^{\nu'}_k)
\end{split}
\end{equation*}
with 
\begin{equation*}\label{Existence_First_Derivative_2}
\begin{split}
\frac{d}{d\e}m_\b(\nu^{\nu'+\e\r}_k)&=\frac{d}{d\e}\frac{\int_{S_k}e_{\tilde\o}e^{\b \langle e_{\tilde\o},M_\b(\nu'+\e\r)\rangle}\a(d\tilde\o)}{\int_{S_k}e^{\b \langle e_{\tilde\o},M_\b(\nu'+\e\r)\rangle}\a(d\tilde\o)}\cr
&=\b(\frac{\int_{S_k}e_{\tilde\o}e^{\b \langle e_{\tilde\o},M_\b(\nu')\rangle}\langle e_{\tilde\o},\frac{d}{d\e}M_\b(\nu'+\e\r)\rangle\a(d\tilde\o)}{\int_{S_k}e^{\b \langle e_{\tilde\o},M_\b(\nu')\rangle}\a(d\tilde\o)}\cr
&\hspace{0,1cm}-\frac{\int_{S_k}e_{\tilde\o}e^{\b \langle e_{\tilde\o},M_\b(\nu')\rangle}\a(d\tilde\o)\int_{S_k}e^{\b \langle e_{\tilde\o},M_\b(\nu')\rangle}\langle e_{\tilde\o},\frac{d}{d\e}M_\b(\nu'+\e\r)\rangle\a(d\tilde\o)}{(\int_{S_k}e^{\b \langle e_{\tilde\o},M_\b(\nu')\rangle}\a(d\tilde\o))^2})\cr
&=:W(k,M_\b(\nu'))\frac{d}{d\e}M_\b(\nu'+\e\r).
\end{split}
\end{equation*}
Here $W(k,M_\b(\nu'))=\b\begin{pmatrix}
 A(k) & B(k)\\
 B(k) & C(k) 
 \end{pmatrix}$ is a $2\times2$ matrix with
\begin{equation*}\label{Mean_Field_Equations_6}
\begin{split}
A(k)&:=\frac{\int_{S_k}\cos^2(\tilde\o) e^{\b\langle e_{\tilde\o},M_\b(\nu')\rangle}\a(d\tilde\o)}{\int_{S_k} e^{\b \langle e_{\tilde\o},M_\b(\nu')\rangle}\a(d\tilde\o)}
-\bigl(\frac{\int_{S_k}\cos(\tilde\o)e^{\b \langle e_{\tilde\o},M_\b(\nu')\rangle}\a(d\tilde\o)}{\int_{S_k}e^{\b \langle e_{\tilde\o},M_\b(\nu')\rangle}\a(d\tilde\o)}\bigr)^2\cr
B(k)&:=
\frac{\int_{S_k}\cos(\tilde\o)\sin(\tilde\o)e^{\b\langle e_{\tilde\o},M_\b(\nu')\rangle}\a(d\tilde\o)}{\int_{S_k} e^{\b \langle e_{\tilde\o},M_\b(\nu')\rangle}\a(d\tilde\o)}\cr
&\hspace{2cm}-\frac{\int_{S_k}\cos(\tilde\o)e^{\b \langle e_{\tilde\o},M_\b(\nu')\rangle}\a(d\tilde\o)\cdot\int_{S_k}\sin(\tilde\o)e^{\b \langle e_{\tilde\o},M_\b(\nu')\rangle}\a(d\tilde\o)}{[\int_{S_k}e^{\b \langle e_{\tilde\o},M_\b(\nu')\rangle}\a(d\tilde\o)]^2}\cr
C(k)&:=\frac{\int_{S_k}\sin^2(\tilde\o) e^{\b\langle e_{\tilde\o},M_\b(\nu')\rangle}\a(d\tilde\o)}{\int_{S_k} e^{\b \langle e_{\tilde\o},M_\b(\nu')\rangle}\a(d\tilde\o)}
-\bigl(\frac{\int_{S_k}\sin(\tilde\o)e^{\b \langle e_{\tilde\o},M_\b(\nu')\rangle}\a(d\tilde\o)}{\int_{S_k}e^{\b \langle e_{\tilde\o},M_\b(\nu')\rangle}\a(d\tilde\o)}\bigr)^2\cr
\end{split}
\end{equation*}
some functions of the covariances. Consequently
\begin{equation}\label{MF_Fixed_Point_Equation}
 \begin{split}
\frac{d}{d\e}M_\b(\nu'+\e\r)&=\sum_l\r(l)[I_{2\times2}-\sum_k\nu'(k)W(k,M_\b(\nu'))]^{-1}m_\b(\nu^{\nu'}_l)
\end{split}
\end{equation}
whenever the matrix inverse exists.
Up to this point all calculations are made for general $\nu'\in \PP(\{1,\dots,q\})$. Evaluating at the equidistribution ($eq$) we find
\begin{equation*}\label{Linearization_RHS_Eq}
 \begin{split}
d F_{eq}(\r)(k)&=\frac{1}{q}(\frac{d}{d\eps}c(k-1,eq + \eps \rho)- \frac{d}{d\eps}c(k,eq + \eps \rho))+\frac{q}{2\pi}(\rho(k-1)- \rho(k)).\cr
\end{split}
\end{equation*}
Since $m_\b(\nu_k^{eq})=\frac{q}{2\pi}\int_{S_k} e_\o\a(d\o)$ we can write
\begin{equation*}
 \begin{split}
\frac{d}{d\e}c(k,eq+\eps\rho)&=\frac{\b q}{2\pi}[\langle e_{\frac{2\pi}{q}k}
,\frac{d}{d\e}M_\b(eq+\e\r)\rangle-\frac{q}{2\pi}\langle\int_{S_k} e_\o\a(d\o),\frac{d}{d\e}M_\b(eq+\e\r)\rangle]\cr
&=\frac{\b q}{2\pi}\langle 
\left (\begin{array}{lllllll}
\cos(\frac{2\pi}{q}k)-\frac{q}{2\pi}\int_{S_k} \cos(\o)\a(d\o)$$
\\ [6 mm]$$\sin(\frac{2\pi}{q}k)-\frac{q}{2\pi}\int_{S_k} \sin(\o)\a(d\o)
\end{array}\right )
,\frac{d}{d\e}M_\b(eq+\e\r)\rangle.\cr
\end{split}
\end{equation*}
Thus for the vector 
\begin{equation*}
 \begin{split}
v(k)&=\left (\begin{array}{lllllll}
\frac{q}{2\pi}(\int_{S_{k}} \cos(\o)\a(d\o)-\int_{S_{k-1}} \cos(\o)\a(d\o))+\int_{S_{k}} \sin(\o)\a(d\o)$$
\\ [6 mm]$$\frac{q}{2\pi}(\int_{S_{k}} \sin(\o)\a(d\o)-\int_{S_{k-1}} \sin(\o)\a(d\o))-\int_{S_{k}} \cos(\o)\a(d\o)
\end{array}\right )
\end{split}
\end{equation*}
which is close to zero for large $q$, we have
\begin{equation*}\label{Linearization_RHS_Eq_2}
 \begin{split}
d F_{eq}(\r)(k)&=\frac{q}{2\pi}(\rho(k-1)- \rho(k))+\frac{\b}{2\pi}\langle v(k),\frac{d}{d\e}M_\b(eq+\e\r)\rangle.\cr
\end{split}
\end{equation*}
\textbf{Remark: } For small $\b$, $d F_{eq}$ is a small pertubation of the rotation matrix $\frac{q}{2\pi}(D-I)$ with $D_{kl}=1_{l=k+1}$. 
Thinking of $\frac{q}{2\pi}$ as a time rescaling, one can consider the linear system of differential equations on probability vectors $x$ of lenght $q$ 
\begin{equation*}\label{rotation_generator2}
\begin{split}
\dot x&=(D-I)x.
 \end{split}
\end{equation*}
Using discrete Fourier transform, it is immediately seen, that this system is attractive towards the equidistribution. 

\bigskip
Let us look at the effect of the pertubation:
\begin{equation*}\label{Local_Magnetisation_3}
 \begin{split}
\frac{d}{d\e}M_\b(eq+\e\r)&=\sum_l\r(l)[I_{2\times2}-\frac{1}{q}\sum_kW(k,0)]^{-1}m_\b(\nu^{eq}_l)
\end{split}
\end{equation*}
where $I_{2\times 2}-\frac{1}{q}\sum_kW(k,0)=\begin{pmatrix}
 1-\frac{\b}{2}(1
 -(\frac{q}{\pi})^2\sin^2(\frac{\pi}{q}))& 0\\
0& 1-\frac{\b}{2}(1
 -(\frac{q}{\pi})^2\sin^2(\frac{\pi}{q}))
 \end{pmatrix}$.
Hence for $\b(1-(\frac{q}{\pi})^2\sin^2(\frac{\pi}{q}))\neq2$ we have
\begin{equation*}\label{Local_Magnetisation_4}
 \begin{split}
\frac{d}{d\e}M_\b(eq+\e\r)&=\frac{q}{2\pi-\b\pi(1
 -(\frac{q}{\pi})^2\sin^2(\frac{\pi}{q}))}\sum_l\r(l)\int_{S_l} e_\o\a(d\o)
\end{split}
\end{equation*}
and thus
\begin{equation*}\label{Linearization_RHS_Eq132}
 \begin{split}
&d F_{eq}(\r)(k)=\frac{q}{2\pi}[(\rho(k-1)- \rho(k))\cr
&\hspace{2.3cm}+\frac{\b}{2\pi-\b\pi(1
 -(\frac{q}{\pi})^2\sin^2(\frac{\pi}{q}))}\sum_l\langle \int_{S_l} e_\o\a(d\o),v(k)\rangle\r(l)]\cr
&=\frac{q}{2\pi}[(\rho(k-1)- \rho(k))+\frac{4\b\sin^2(\frac{\pi}{q})}{2\pi-\b\pi(1
 -(\frac{q}{\pi})^2\sin^2(\frac{\pi}{q}))}\sum_l\sin(\frac{2\pi}{q}(k-l))\r(l)\cr
&\hspace{0.5cm}+\frac{2\b q\sin^2(\frac{\pi}{q})}{2\pi^2-\b\pi^2(1
 -(\frac{q}{\pi})^2\sin^2(\frac{\pi}{q}))}\sum_l(\cos(\frac{2\pi}{q}(k-l))-\cos(\frac{2\pi}{q}(k-l-1)))\r(l)].\cr
\end{split}
\end{equation*}
In matrix notation this is
\begin{equation*}\label{Linearization_RHS_Eq_10}
 \begin{split}
\tilde M_\b(i,j):=\frac{q}{2\pi}(\d_{j=(i-1)}-\d_{j=i}&+c_1\sin(\frac{2\pi}{q}(i-j))\cr
&+c_2[\cos(\frac{2\pi}{q}(i-j))-\cos(\frac{2\pi}{q}(i-j-1)])\cr
\end{split}
\end{equation*}
where the property $\sum_j \tilde M_\b(i,j)=0$ for all $i\in\{1,\dots,q\}$ reflects conservation of mass. 
\begin{lem}\label{Eigenvalues_Lemma}
The eigenvalues of $\tilde M$ are given by
\begin{equation*}\label{Eigenvalue_2}
 \begin{split}
&\l_1=\frac{q}{2\pi}([(c_2\frac{q}{2}-1)(1-\cos(\frac{2\pi}{q}))]+i[(c_2\frac{q}{2}-1)\sin(\frac{2\pi}{q})-c_1\frac{q}{2}])\cr
&\l_j=\frac{q}{2\pi}([\cos(\frac{2\pi}{q}j)-1]\pm i\sin(\frac{2\pi}{q}j))\hspace{1cm}\text{for }j\in\{2,\dots,q-2\}\cr
&\l_{q-1}=\frac{q}{2\pi}([(c_2\frac{q}{2}-1)(1-\cos(\frac{2\pi}{q}))]-i[(c_2\frac{q}{2}-1)\sin(\frac{2\pi}{q})-c_1\frac{q}{2}])\cr
&\l_q=0.\cr
\end{split}
\end{equation*}
where $c_1=\frac{4\b\pi\sin^2(\frac{\pi}{q})}{2\pi^2-\b\pi^2+\b q^2\sin^2(\frac{\pi}{q})}$ and $c_2=\frac{2\b q \sin^2(\frac{\pi}{q})}{2\pi^2-\b\pi^2+\b q^2\sin^2(\frac{\pi}{q})}$.
\end{lem}
\textbf{Proof of Lemma \ref{Eigenvalues_Lemma}: }
Since $\tilde M$ is rotation invariant, we can employ discrete Fourier transformation to calculate the eigenvalues and eigenvectors of $\tilde M$. The $k$-th eigenvector is given by
\begin{equation*}\label{Eigenvector}
 \begin{split}
u_k=\frac{1}{\sqrt{q}}\exp(i\frac{2\pi}{q}kl)_{l\in\{1,\dots,q\}}
\end{split}
\end{equation*}
and with $\hat M_\b(n):=\tilde M_\b(i,(i+n)\text{mod}(q))$ the $k$-th eigenvalue reads
\begin{equation*}\label{Eigenvalue}
 \begin{split}
\l_k=\sum_{n=0}^{q-1}\hat M_\b(n)\exp(-i\frac{2\pi}{q}kn).
\end{split}
\end{equation*}
Calculating separately for the summands in $\tilde M$, the result follows.
$\Cox$

\bigskip
Notice, the eigenvalues always come in conjugated pairs. The eigenvectors have zero weight (except for the one belonging to the zero eigenvalue.)

\bigskip
\textbf{Proof of Theorem \ref{Attract_Equi}: } The eigenspaces for the eigenvalues $\l_2,\dots,\l_{q-2}$ have negative real part and therefore belong to the attractive manifold of the equidistribution. 
The eigenspaces for the perturbated eigenvalues $\l_1,\l_{q-1}$ form a locally non-attractive manifold if 
\begin{equation*}\label{Eigenvalue_3}
 \begin{split}
&\frac{2\pi}{q}\text{Re}(\l_{1,q-1})
=(\cos(\frac{2\pi}{q})-1)(1-\frac{q}{2}c_2)>0
\end{split}
\end{equation*}
or equivalently if $\frac{q}{2}c_2=\frac{\b q^2 \sin^2(\frac{\pi}{q})}{2\pi^2-\b\pi^2+\b q^2\sin^2(\frac{\pi}{q})}>1$. Since we assume $\b>2$ this is again equivalent to $2>\b(1-(\frac{q}{\pi})^2\sin^2(\frac{\pi}{q}))$. 
Now if $2>\b(1-\frac{q}{2\pi}\sin(\frac{2\pi}{q}))$ (in other words \eqref{Non_Gibbs_Criterion} fails and we are in the relevant parameter regimes), then $2>\b(1-(\frac{q}{\pi})^2\sin^2(\frac{\pi}{q}))$. 
But this is true since $\b(1-\frac{q}{2\pi}\sin(\frac{2\pi}{q}))>\b(1-(\frac{q}{\pi})^2\sin^2(\frac{\pi}{q}))$ is equivalent to $\cos(\frac{\pi}{q})<\frac{q}{\pi}\sin(\frac{\pi}{q})$  and $\frac{q}{\pi}\sin(\frac{\pi}{q})=\cos(\xi)$ for some $\xi\in[0,\frac{\pi}{q}]$.
$\Cox$

\bigskip
An illustration is given in figure~\ref{Plot_Parameter_Regimes}. Notice, $\lim_{q\to\infty}\frac{q}{2}c_2=\frac{\b}{2}$ and hence all real parts go to zero as they should.
We would like to point out, that in \cite{GPPP12} although a rotation dynamics on the continuous system driven by Brownian motion is considered, similar attractivity conditions appear. In particular, in the low temperature regime, the periodic orbit attracts every measure, except the equidistribution and whatever is attracted to it. The attractive manifold for the equidistribution is also given by a continuous version of    
\begin{equation*}\label{Attractiv_Manifold}
 \begin{split}
\{\nu'\in \PP(\{1,\dots,q\})|\sum_{k=1}^q\nu'(k)\exp(i\frac{2\pi}{q}k)=0\text{ and }\sum_{k=1}^q\nu'(k)\exp(-i\frac{2\pi}{q}k)=0\}.
\end{split}
\end{equation*}

\end{document}